\documentstyle[11 pt, amscd,amstex, epsfig, color]{article}
\input xypic
\xyoption{all}

\def\M{{\cal M}}
\def\CM{\overline{\cal M}}

\def\H{{\cal H}}
\def\F{\overline{\cal F}}

\def\L{{\cal L}}
\def\J{{\cal J}}

\def\a{\alpha}
\def\o{\omega}
\def\l{\lambda}

\def\b{\beta}
\def\dbar{\overline{\partial}}

\def\non{\noindent}
\def\pf{\non {\bf Proof. }}
\def\QED{\nopagebreak \hskip .1in { $\Box$ }\penalty10000 %
\hskip\parfillskip \par  }

\newtheorem{theorem}{Theorem}[section]
\newtheorem{prop}[theorem]{Proposition}
\newtheorem{lemma}[theorem]{Lemma}

\newtheorem{defn}[theorem]{Definition}

\newtheorem{rem}[theorem]{Remark}

\def\dga{\M^{V,\delta}_{g,k}(\,d\,)}
\def\dgb{\CM^{\,\delta}_{g,k}(\,\l,d\,)}
\def\dgc{\bigcup\ \big(\,\CM_{g_{1},k_{1}+1}\times\CM_{g_{2},k_{2}+1}\,\big)
          \times\big(\,E(n)^{k_{1}}\times E(0)^{k_{2}}\,\big)}
\def\dgcc{\big(\,\CM\times E(n)\,\big)_{g,k}}
\def\dgd{\CM_{g,k}\times Z_{\l}^{k}}
\def\dge{\CM_{g,k}\times Z_{0}^{k}}

\begin{document}

%%%%%%%%%%%%%%%%%%%%%%%%%%%%%%%%%%%%%%%
%%%%%%%%%%%%%%%%%%%%%%%%%%%%%%%%%%%%%%%
\date{\empty}
\title{\bf  Counting Curves in Elliptic Surfaces by Symplectic Methods \\
      \vskip.2in}
\author{Junho Lee\\
University of Minnesota\\
Minneapolis, MN, 55454}

\addtocounter{section}{0}

%%%%%%%%%%%%%%%%%%%%%%%%%%%%%%%%%%%%%%%

\maketitle

\vskip.15in
%%%%%%%%%%%%%%%%%%%%%%%%%%%%%%%%%%%%%%%

\begin{abstract}
We explicitly compute family GW invariants of elliptic surfaces for primitive classes.  That involves   establishing a
TRR formula and a symplectic sum formula for elliptic surfaces and  then determining the GW invariants using an argument from \cite{ip3}.  In particular, as  in \cite{bl1},  these calculations  also
confirm the well-known Yau-Zaslow Conjecture \cite{yz}
for primitive classes in  $K3$ surfaces.
\end{abstract}

\vskip.4in

In \cite{l} we introduced ``family GW invariants'' for K\"{a}hler
surfaces with $p_g>0$. Since these invariants are defined by using
{\em non-compact} family of almost K\"{a}hler structures, we can
easily extend several existing techniques for calculating GW
invariants to the family GW invariants. In particular, the `TRR
formula' applies to the family invariants, and at least some
special cases of the symplectic sum formula \cite{ip3} apply, with
appropriate minor modifications to the formula.  Those  formulas
enable us  to enumerate the curves in the elliptic surfaces $E(n)$
for the class $A$= section plus multiples of the fiber.

\bigskip

\begin{theorem}
Let $E(n)\to {\Bbb P}^{1}$  be a standard elliptic surface with a
section of self-intersection $-n$. Denote by $s$ and $f$ the
homology class of the section and the fiber.  Then the genus $g$
family GW invariants for the classes $s+df$ are given by the
generating function
\begin{equation}\label{0.2}
\sum\limits_{d\geq 0}\,
           GW_{s+df,g}^{\H}(E(n))(\,pt^{g}\,)\,t^{d}\
           =\ \left(\,tG^{\prime}(t)\,\right)^{g}\,
            \prod\limits_{d\geq 1}
              \left(\,1-t^{d}\, \right)^{-12n}
\end{equation}
where $G(t) = \sum\limits_{d\geq 1}\sigma(d)\,t^{d}$\ \  and\ \
      $\sigma(d)=\sum_{k|d}k$\,.
\label{thm02}
\end{theorem}

\medskip

Bryan and Leung (\cite{bl1},\cite{bl2}) defined family invariants
for K3 and Abelian surfaces by using the Twistor family. They used
algebraic methods to show (\ref{0.2}) for GW invariants of the
rational elliptic surface $E(1)$ and for family invariants of
$E(2)=K3$ surfaces. For K3, that confirms  the famous Yau-Zaslow
Conjecture \cite{yz} for those cases when the homology class $A$
is primitive. They also pointed out that one can define family
invariants of $E(n)$ for $n\geq 3$ using compact family of complex
structures induced from the fiber sum, and then use the algebraic
methods  of \cite{bl1} to show that those invariants also satisfy
(\ref{0.2}) \cite{bl4}, see also section 5 of \cite{bl3}.

\smallskip
On the other hand, Ionel and Parker used analytic methods to
compute the GW invariants of $E(1)$ \cite{ip3}. They related TRR
formula and their sum formula for the relative invariants to
obtain a quasi-modular form as in (\ref{0.2}). We follow the same
argument --- relating TRR formula and sum formula --- to show
Theorem~\ref{thm02}. This theorem also confirm the Yau-Zaslow
Conjecture for primitive classes, since our invariants of K3
surfaces are equivalent to the invariants define by Bryan and
Leung\,(cf. Theorem 4.3 of \cite{l}).
\bigskip

 For $E(1)$ and $E(2)$, the  invariants are known to be enumerative,
 that is, formula  (\ref{0.2}) actually counts
 (irreducible) holomorphic curves in the primitive classes
 for generic complex structures on those surfaces \cite{bl1}.
 At the moment, it is not clear whether, or in what sense,
 that is true for the $E(n)$ with
 $n\geq 3$\,(cf. Remark 5.12 of \cite{bl3}).

\vskip 0.6 cm

The construction of family invariants
for K\"{a}hler surfaces is briefly described in Section 1.
We give  an overview of the proof of
Theorem~\ref{thm02} in Section 2.
This  argument is an  extension of the elegant argument
used by Ionel and Parker to
compute the  GW invariants of $E(1)$ \cite{ip3}.  It involves
computing the generating function for the
  invariants in two ways, first using the so-called TRR formula, and
second using a symplectic sum formula
as in \cite{ip3}.
Roughly, the only  modification  needed is a shift
in the dimension counts.
The extended TRR formula is proved in Section 3 and
the sum formula are
established in the last 5 sections.

\smallskip

Section 4 gives an alternative definition of the family  invariants
for $E(n)$ based on the idea of  perturbing  the
$J_{\a}$-holomorphic map equations
as in \cite{rt1,rt2}.
This alternative definition is better suited  to adapt the
analytic arguments in \cite{ip2,ip3}
to a family version of sum formula.  The proof of the sum formula
begins by studying holomorphic maps
into a degeneration of $E(n)$.   Because $E(n)$ is a K\"{a}hler
surface we are able to
degenerate within a holomorphic family, rather than the symplectic
family used in \cite{ip3}.

\smallskip

The degeneration family $Z$ is described in Section 5.  It is a
family $\lambda : Z\to D^{2}$ whose
  fiber $Z_\l$ at $\l\neq 0$ is a copy of $E(n)$ and whose central
fiber is   a union of $E(n)$ with
$E(0)=T^{2}\times S^{2}$ along a fixed elliptic fiber $V$.
%Gromov Convergence Theorem \cite{is,p,pw} and
%compactness of family moduli space shown in Section 6 imply that
As $\l\to 0$ maps into $Z_\l$ converge to maps into $Z_0$, and by
bumping  $\alpha$ to zero along the fiber $V$ we can ensure that
the limits satisfy a simple matching condition along $V$ (there is
a single matching condition for the classes $A$ that we consider).
Section 6 shows this splitting argument.

\smallskip
Conversely, if a map into $Z_{0}$ satisfies the matching condition
then it can be smoothed to produce a map into $Z_{\lambda}$ for
small $\l$. That smoothing is the Gluing Theorem in \cite{ip3},
which relate family invariants of $E(n)$ with relative invariants
of $E(n)$ and $E(0)$ relative to $V$. We define a family version
of relative invariants of $E(n)$ in Section 7. Using the Gluing
Theorem, we prove the required sum formulas for the  family
invariants of $E(n)$ in Section 8.

\vskip 0.6cm \noindent {\bf Acknowledgements : } I would like to
thank most sincerely my advisor Prof. Thomas Parker for his
guidance and helpful discussions. Without his help this paper
would not have been possible. I also wish to thank Prof. Naichung
Leung and Prof. Eleny Ionel for useful discussions.

\vskip 1cm

%%%%%%%%%%%%%%%%%%%%%%%%%%%%%%%%%%%%%%%%%%%%%%%%%%%%%%%%%%%%
%%%%%%%%%%%%%%%%%  Section 1%%%%%%%%%%%%%%%%%%%%%%%%%%%%%%
%%%%%%%%%%%%%%%%%%%%%%%%%%%%%%%%%%%%%%%%%%%%%%%%%%%%%%%%%%%%%%
\setcounter{equation}{0}
\section{Family Invariants for K\"{a}hler surfaces}
\label{section2}
\bigskip

Let $X$ be a closed complex surface with  K\"{a}hler structure
$(\o,J,h)$.  In this section we briefly describe the family Gromov-Witten invariants  associated to $(X,J)$ which were defined in  \cite{l}.  First set
  \begin{equation*}\label{para}
    \H\ =\ \{\ \a+\overline{\a}\ |\ \a\,\in\,H^{2,0}(X)\ \}.
  \end{equation*}
This  is a  $2\,p_{g}$-dimensional space of harmonic forms which are
  {\em $J$-anti-invariant}, that is,
$\a(Ju,Jv)= -\a(u,v)$.
Each  $\a\in\H$ defines an endormorphism
$K_{\a}$ of $TX$ by the equation
  \begin{equation*}
    h(u,K_{\a}v )\  =\  \a(u,v).
  \end{equation*}
One can check  that, for each $\a\in\H$,  $Id+JK_{\a}$ is invertible, so  defines an almost complex structure
  \begin{equation*}
    J_{\a}\ =\ \left(\, Id+JK_{\a}\,\right)^{-1}\,J\,
               \left(\, Id+JK_{\a}\,\right).
  \end{equation*}

\medskip
Let $\F=\F_{g, k, A}$ the space of all stable maps
$f:(C,j)\to X$ of genus $g$ with $k$ marked points which
represent homology class $A$.
For each such map, collapsing unstable components of the domain
determines a point
in the Deligne-Mumford space $\CM_{g,k}$ and
evaluation of marked points determines a point in $X^{n}$.
Thus we have a map
  \begin{equation}\label{st-ev}
     \F\, @>\  st\times ev\  >> \,\CM_{g,k}\times X^{k}
  \end{equation}
where $st$ and $ev$ denote stabilization map and evaluation maps,
respectively.
On the other hand, there is a {\em generalized orbifold bundle}
$E$ over $\F\times \H$ whose
fiber over $(\,f,j,\a\,)$ is
$\Omega^{0,1}_{jJ_{\a}}(f^{*}TX)$.
This bundle has a section $\Phi$ defined by
  \begin{equation}\label{defnSection}
     \Phi(f,j,\a)=df+J_{\a}\,df\,j\,.
  \end{equation}
By definition, the right-hand side of (\ref{defnSection})
vanishes for $J_{\a}$-holomorphic maps.
Thus $\Phi^{-1}(0)$ is the moduli space of
$J_{\a}$-holomorphic maps which we denote by
  \begin{equation*}
    \CM_{g,k}^{J,\H}(\,X,A\,).
  \end{equation*}
It is, unfortunately, not always compact.  When it is compact, it gives rise to family Gromov-Witten invariants in the usual way (cf. \cite{l}).

\medskip
\begin{prop}{\rm (\cite{lt})}\label{EulerClass}
Suppose the moduli space $\Phi^{-1}(0)$ is compact.
Then the bundle $E$ has a rational homology ``Euler class''  $[\,\CM_{g,k}^{J,\H}(X,A)\,]^{\rm vir} \in H_{2r}( \F;{\Bbb Q})$ for
  \begin{equation*}
     r\ =\ c_{1}(X)[A]\ +\ g-1\ + \ k\ +\ p_{g}.
  \end{equation*}
\end{prop}

%The proof of Proposition~\ref{EulerClass} follows from
%Proposition 2.2 and Theorem 1.2 of \cite{lt}.

\medskip
\begin{defn}\label{D:fgw-inv't}
Whenever the moduli space $\CM_{g,k}^{J,\H}(\,X,A\,)$ is compact,
we define the family GW invariants of $(X,J)$ to be the map
  \begin{equation*}
     GW^{J,\H}_{g,k}(X,A)\,:\,
        H^{*}(\,\CM_{g,k};{\Bbb Q}\,)\times
        [\,H^{*}(\,X;{\Bbb Q}\,)\,]^{k}\, \to \,{\Bbb Q}
  \end{equation*}
defined on  $\beta\in H^{*}(\,\CM_{g,k};{\Bbb Q}\,)$ and
      $\a\in H^{*}(\,X^{k};{\Bbb Q}\,)$ by
  \begin{equation*}
     GW^{J,\H}_{g,k}(X,A)(\,\beta;\a\,)\ =\
     [\,\CM_{g,k}^{J,\H}(X,A)\,]^{\rm vir}\,\cap\,
     \left(\, st^{*}(\beta)\,\cup\,ev^{*}(\a)
     \,\right).
  \end{equation*}
\end{defn}

\bigskip

This paper will focus on the case where $X$ is
a standard elliptic surface $E(n)$ with a section class.
Note that the elliptic surfaces $E(2)$ are K3 surfaces.
Denote by  $s$ and $f$ the homology class of the section
and the fiber of $E(n)$.
Since $c_1(E(n))=(2-n)f$ and $p_g=n-1$,  we have
\begin{equation}
\label{1.r}
\mbox{dim } \CM_{g,k}^{J,\H}(\,E(n),A\,)\ =\ 2(g+k).
\end{equation}

\medskip
\begin{prop}{\rm (\,\cite{l}\,)}\label{P:Elliptic}
Let $(X,J)$ be an elliptic surface $E(n)$ and
$A=s+df$, where $d$ is an integer.
  \begin{enumerate}
    \item[(a)] The moduli space $\CM^{J,\H}_{g,k}(X,A)$
               is compact and hence
               the invariants $GW^{J,\H}_{g,k}(X,A)$ are well-defined.
    \item[(b)] The invariants $GW^{J,\H}_{g,k}(X,A)$
               depend only on the deformation class
               of $(X,J)$.
    \item[(c)] For K3 surfaces (i.e. n=2) the
               $GW^{J,\H}_{g,k}(X,A)$ are same as the invariants defined by
               Bryan and Leung in {\rm \cite{bl1}}.
  \end{enumerate}
\end{prop}

Thus for elliptic surfaces the family of $J_{\a}$-holomorphic maps parameterized by the family $\H$ gives rise to well-defined invariants, which we will denote variously as
  \begin{equation*}
    GW^{\H}_{g,k}(E(n),A)\, ,\ \
    GW^{\H}_{A,g}(E(n))\, ,
    \mbox{\ \ or\ simply\ \ }
    GW^{\H}_{A,g}.
  \end{equation*}
The goal of this paper is to {\em calculate} these family GW invariants.

\medskip

The family invariants have a property analogous to
the composition law of ordinary GW invariants.
Consider a node $p$ of a stable curve $C$
in the Deligne-Mumford space $\CM_{g,k}$.
When the node is separating, the normalization of $C$ has two components.  The
genus and the number of marked points decompose as
$g=g_{1}+g_{2}$ and $k=k_{1}+k_{2}$ and
there is a natural map
  \begin{equation}\label{gluing-1}
     \sigma : \CM_{g_{1},k_{1}+1}\times\CM_{g_{2},k_{2}+1}\to
              \CM_{g,k}.
  \end{equation}
defined by identifying $(k_{1}+1)$-th marked points
of the first component to the first marked point of the second component.
We denote by $PD(\sigma)$ the Poincar\'{e} dual of the image of
this map $\sigma$.
For non-separating node, there is another natural map
  \begin{equation*}\label{gluing-2}
    \theta:\CM_{g-1,k+2}\to \CM_{g,k}
  \end{equation*}
defined by identifying the last two marked points.
We also write $PD(\theta)$ for the Poincar\'{e} dual of the image of
this map $\theta$.

\begin{prop}{\rm (\cite{l})} \label{composition}
Let $\{H_{\gamma}\}$ be any basis of $H^{*}(X;{\Bbb Z})$ and
$\{H^{\gamma}\}$ be its dual basis and suppose that
$GW^{J,\H}_{g,k}(X,A)$ is defined.
  \begin{enumerate}
    \item[{\rm(a)}] Given any decomposition $A=A_{1}+A_{2}$,\
                    $g=g_{1}+g_{2}$, and $k=k_{1}+k_{2}$,
                    if the moduli space
                    $\CM_{g_{1},k_{1}}^{J,\H}(\,X,A_{1}\,)$
                    is compact, then
                    \begin{align*}
                      & GW^{J,{\cal H}}_{A,g}(X)(PD(\sigma);
                                        \a_{1},\cdots,\a_{k}) \\
                      & =
                      \sum_{A=A_{1}+A_{2}}\sum_{\gamma}
                      GW^{J,\H}_{A_{1},g_{1}}(X)
                       ( \a_{1},\cdots,\a_{k_{1}},H_{\gamma})\,
                      GW_{A_{2},g_{2}}(X)
                       (H^{\gamma},\alpha_{k_{1}+1},
                               \cdots,\alpha_{k})
             \end{align*}
                   where $GW_{A_{2},g_{2}}(X)$
                   denotes the ordinary GW invariant.
     \item[(b)] $  GW^{J,{\H}}_{A,g}(X)(PD(\theta);
                                        \alpha_{1},\cdots,\alpha_{k}) =
                   \sum_{\gamma}GW^{J,{\H}}_{A,g-1}(X)
                    (\a_{1},\cdots,\a_{k},H_{\gamma},H^{\gamma}).$
\end{enumerate}
\end{prop}

\vskip 1cm
%%%%%%%%%%%%%%%%%%%%%%%%%%%%%%%%%%%%%%%%%%%%%%%%%%%%%%%%%%%%%%%%
%%%%%%%%%%%%%% Section 2 %%%%%%%%%%%%%%%%%%%%%%%%%%%%%%%%%%%%%%%
%%%%%%%%%%%%%%%%%%%%%%%%%%%%%%%%%%%%%%%%%%%%%%%%%%%%%%%%%%%%%%%%%
\setcounter{equation}{0}
\section{The Invariants of $E(n)$ --- Outline}
\label{section2}
\bigskip

By Proposition~\ref{P:Elliptic},
the family GW invariants of $E(n)$ for the class $s+df$
are unchanged under deformations of K\"{a}hler structure.
Since the moduli space with genus $g$ and no marked points has  dimension $2g$, we get numerical invariants by imposing $g$ point constraints on the moduli space $\CM_{g,g}^{\H}(\,E(n), s+df)$ --- those are the numbers we aim to calculate. For convenience we assemble them in  the generating function
  \begin{equation}
  \label{6.defF}
    F_{g}(t)=\sum\limits_{d\geq 0}\,
    GW_{s+df,g}^{\H}(E(n))(\,pt^{g}\,)\,t^{d}.
  \end{equation}
In  this and the following four sections we will derive the formula for $F_{g}(t)$ stated in Theorem
\ref{0.2}.  Thus our aim
it to prove:

%%%% Main Formula label-Main %%%%%%5
%
%
%

\begin{prop}  \label{P:C-Main}
For $n\geq 1$,
\begin{equation}\label{E:Main}
F_{g}(t)\ =\  \left(\,tG^{\prime}(t)\,\right)^{g}\,
              \prod\limits_{d\geq 0}
              \left(\,1-t^{d} \,\right)^{-12n}
\end{equation}
\end{prop}

\medskip

This section shows how Proposition \ref{P:C-Main} follows
from three formulas,
equations (\ref{E:TRR}),\ (\ref{E:Sum1})\, and (\ref{E:Sum2}) below,
that are proved in later sections.
Our proof parallels the proof of Ionel and Parker
for GW invariants of $E(1)$ \cite{ip3}.

\bigskip

Here is the outline the proof of (\ref{E:Main}).
Consider the `descendent'
$\tau(f^{*})= \psi_{1}\,\cup\, ev^{*}(f^{*})$\  where
$\psi_{1}$ denotes the first Chern class
of the line bundle\ \
$L\to  \overline{{\cal M}}_{1,1}^{\H}(\,E(n),s+df\,)$
whose geometric fiber over
$\left(\,f,(C;x),\a\,\right)$ is $T^{*}_{x}C$.
We define  the generating function for a
genus 1 invariant with the descendent constraint,
namely
\begin{equation}
H(t)    \ =\  \sum\limits_{d\geq 0}
         GW_{s+df,1}^{\H}(E(n))
                  (\,\tau(f^{*})\,)\,t^{d}.
\label{defH(t)}
\end{equation}

We can compute $H(t)$ in two different ways. In section 3,
we show how to
combine the composition law together with the TRR for genus 1
to obtain the formula
%%%%%%%%  TRR   label-TRR   %%%%%%%%
%
  \begin{equation}\label{E:TRR}
     H(t)\  =\  \frac{1}{12}\,t\,F_{0}^{\prime}(t)\  -\
                \frac{1}{12}\,F_{0}(t)\ +\  (2-n)\,F_{0}(t)\,G(t)
  \end{equation}
Then, from section 4 to 8 we establish a family version of
the sum formulas to show
%%%%%%%%  SumFormular label-Sum %%%%%
%
  \begin{align}
     H(t)\ & = \ -\frac{1}{12}\,F_{0}(t)\  +\  2\,F_{0}(t)\,G(t)
     \label{E:Sum1} \\
     F_{g}(t)\ &=\  F_{g-1}(t)\,t\,G^{\prime}(t)
     \label{E:Sum2}
  \end{align}
(see Proposition \ref{Lastprop}).
Equations (\ref{E:TRR}) and (\ref{E:Sum1})
give rise to the ODE
  \begin{equation}\label{ODE}
     t\,F_{0}^{\prime}(t)\  =\  12\,n\,G(t)\,F_{0}(t)
  \end{equation}
and we show in Proposition \ref{F(0)=1} that the initial condition
is  $F_{0}(0)=1$.
It is well-known  that the solution of this ODE
is given by
\begin{equation*}
F_{0}(t) = \prod\limits_{d\geq 0}
            \left(\, 1-t^{d} \,\right)^{-12n}.
\end{equation*}
Now, (\ref{E:Sum2}) gives (\ref{E:Main}) by induction.
That completes the proof of Proposition \ref{P:C-Main} and hence of
the main Theorem
\ref{0.2} of the introduction.  The heart of the matter, then, is  to establish formulas  (\ref{E:TRR}),\ (\ref{E:Sum1})\, and (\ref{E:Sum2}).

\vskip 1cm

%%%%%%%%%%%%%%%%%%%%%%%%%%%%%%%%%%%%%%%%%%%%%%%%%%%%%%%%%%%%
%%%%%%%%%%%%%%%%%  Section 3 %%%%%%%%%%%%%%%%%%%%%%%%%%%%%%
%%%%%%%%%%%%%%%%%%%%%%%%%%%%%%%%%%%%%%%%%%%%%%%%%%%%%%%%%%%%%%

\setcounter{equation}{0}
\section{The Topological Recursion Relation (TRR)}
\label{section3}
\bigskip

This section shows the TRR formula (\ref{E:TRR}).
Following \cite{ac},
we denote by ${\cal M}(G)$
the moduli space of all  genus $g$ stable
curves with $k$ marked points
whose dual graph is $G$.
We also denote by $\delta_{G}$
the orbifold fundamental class of
$\overline{{\cal M}(G)}$, that is,
the fundamental class divided by the order
of the automorphisms of a general element of  ${\cal M}(G)$.
Graphs with one edge correspond to degree two classes.
There are two types of such graphs,
one of which is the graph $G_{irr}$ with one vertex of genus $g-1$.

\smallskip
The following is the well-known genus 1 topological recursion relation:
  \begin{equation}\label{TRR}
    \phi_{1}\ (\,= c_{1}(\L)\,)\ =\
    \frac{1}{12}\,\delta_{G_{irr}}
    \mbox{\ \ \ \ in\ \ \ \ }
    H^{2}(\,\CM_{1,1};{\Bbb Q}\,)
  \end{equation}
where the line bundle $\L\to \CM_{1,1}$ has the geometric fiber
$T_{x}^{*}C$ at the point $(C,x)$.

\medskip

%%%%%%%%%    TRR Formula label-TRR  %%%%%%%%%%%
%
\begin{prop}\label{P:TRR}  The generating function
   (\ref{defH(t)}) satisfies
   \begin{equation*}
       H(t)\  =\  \frac{1}{12}\,t\,F^{\prime}_{0}(t)
                  \  -\  \frac{1}{12}\,F_{0}(t)\
              +\  (2-n)\,F_{0}(t)\,G(t).
   \end{equation*}
\end{prop}

\pf
Let $\CM_{0,2}$ be the space of prestable curves of
genus 0 with two marked points \cite{g} and
$\sigma\,:\,\CM_{0,2}\,\times\, \CM_{1,1}\to \CM_{1,1}$ be
the gluing map as in (\ref{gluing-1}).
For any decomposition of $s+df=A_{1}+A_{2}$, we denote by
   \begin{equation}\label{bdry}
      \CM(\,\sigma(A_{1},A_{2})\,)\ \subset\
      \CM_{1,1}^{\H}(s+df)
   \end{equation}
the set of all $(f,C,\a)$ in $\CM_{1,1}^{\H}(A)$
such that
(i)   $C = \sigma(\,C_{1},C_{2}\,)$ for some $C_{1}\in \CM_{0,2}$
      and $C_{2}\in\CM_{1,1}$,
(ii)  the restriction of $f$ to $C_{1}$
      represents $A_{1}$, and
(iii) the restriction of $f$ to $C_{2}$ represents $A_{2}$. By the
machinery of Li and Tian \cite{lt}, there is a virtual fundamental
class
   \begin{equation}\label{bvfc}
      [\,\CM(\,\sigma(A_{1},A_{2})\,)\,]^{\rm vir}
   \end{equation}
associated with (\ref{bdry}) such that
the coefficients
$GW^{\H}_{s+df,1}(\,\tau(f^{*})\,)$
of \,$H(t)$\, are
  \begin{align}\label{A:TRR}
     [\,\CM_{1,1}^{\H}(\,s+df\,)\,]^{\rm vir}\ \cap\  \tau(f^{*})
     \ =&\
     [\,\CM_{1,1}^{\H}(\,s+df\,)\,]^{\rm vir}\ \cap\
     (\,st^{*}\phi_{1}\,\cup\, ev^{*}(f^{*})\,)\notag \\
      &+\
     \sum\,
     [\,\CM(\,\sigma(A_{1},A_{2})\,)\,]^{\rm vir}\ \cap\
     ev^{*}(f^{*})
  \end{align}
where the sum is over all
decompositions $s+df=A_{1}+A_{2}$.

\medskip
Let $\{\,H^{\gamma}\,\}$ and $\{\,H_{\gamma}\,\}$
be bases of $H^{*}(\,E(n)\,)$ dual by the intersection form.
We have
  \begin{align}\label{SC}
     [\,\CM_{1,1}^{\H}(\,s+df\,)\,]^{\rm vir}\ \cap\
     (\,st^{*}\phi_{1}\,\cup\, ev^{*}(f^{*})\,)\
     &=\ \frac{1}{12}\,
     GW^{\H}_{s+df,1}(\,\delta_{G_{irr}};f^{*}\,)\  \notag\\
     &=\ \frac{1}{24}\,\sum_{\gamma}\,
     GW^{\H}_{s+df,0}(\,f^{*},H^{\gamma},H_{\gamma}\,)\ \notag\\
     &=\ \frac{2d-n}{24}\,
     GW^{\H}_{s+df,0}
  \end{align}
where the first equality follows from (\ref{TRR}),
      the second follows from Proposition~\ref{composition}\,b and
      $|\,{\rm Aut}(G_{irr})\,|=2$, and the last follows from
      $\sum_{\gamma}(H^{\gamma}\cdot A)(H_{\gamma}\cdot A)=A^{2}$.

\medskip

On the other hand, by Theorem 2.4 of \cite{l}, every $(f,\a)$ in
$\CM_{1,1}^{\H}(\,s+df\,)$\ has $\a=0$, i.e. $f$ is truly
holomorphic. The inclusion (\ref{bdry}) thus means that the only
possible decompositions of $s+df$ with nontrivial virtual class
(\ref{bvfc}) are  $s+d_{1}f$ and $d_{2}f$ with $d_{1}+d_{2}=d$,
and $d_{1},d_{2}\geq 0$.

 Proposition~\ref{composition}\,a and routine dimension counts then imply that

\begin{equation}
\label{FC}
    \sum\,[\,\CM(\,\sigma(A_{1},A_{2})\,)\,]^{\rm vir}\, \cap\,
    ev^{*}(f^{*})   \ = \
    \sum_{d_{1}+d_{2}}\sum_{\gamma}
    GW^{\H}_{s+d_{1}f}(f^{*},H^{\gamma})\,GW_{d_{2}f,1}(H_{\gamma}).
 \end{equation}
This can be further simplified by separating the $d_2=0$ term and
simplifying using the facts (a)   $\sum_{\gamma}(H^{\gamma}\cdot
A)(H_{\gamma}\cdot B)=A\cdot B$, (b)
$d_{2}\,GW_{d_{2}f,1}=(2-n)\sigma(d_{2})$\,(see \cite{ip1}), and
(c) $GW_{0,1}(H_{\gamma})=\frac{1}{24}(K\cdot H_{\gamma})$
 where $K=(n-2)f$ is the canonical class\,(see 1.4.1
 Proposition of \cite{km1}).  The righthand side of
(\ref{FC}) then becomes
 \begin{equation}
 \label{FC2}
  (2-n)\,   \sum_{1\leq d_{2}\leq d}GW^{\H}_{s+d_{1}f,0}\,\sigma(d_{2})
    \ +\ \frac{n-2}{24}\,GW^{\H}_{s+df,0}
 \end{equation}
\medskip
The  proof now follows from (\ref{A:TRR}),
(\ref{SC}), (\ref{FC2}) and the definitions of $F_{0}(t)$ and $H(t)$.
\QED

\vskip 1cm

%%%%%%%%%%%%%%%%%%%%%%%%%%%%%%%%%%%%%%%%%%%%%%%%%%%%%%%%%%%%
%%%%%%%%%%%%%%%%%  Section 4 %%%%%%%%%%%%%%%%%%%%%%%%%%%%%%
%%%%%%%%%%%%%%%%%%%%%%%%%%%%%%%%%%%%%%%%%%%%%%%%%%%%%%%%%%%%%%

\setcounter{equation}{0}
\section{Ruan-Tian Invariants of $E(n)$}
\label{section4}
\bigskip

Instead of constructing virtual fundamental class directly from
the moduli space of stable $J$-holomorphic maps,
Ruan and Tian \cite{rt1,rt2}
perturbed $J$-holomorphic equation to  $\dbar f=\nu$ where the
inhomogeneous term $\nu$ can be
chosen generically.
For generic $(J,\nu)$, the moduli space
of stable $(J,\nu)$-holomorphic maps is then a
compact smooth orbifold with all
lower strata having codimension at least two.
Ruan and Tian defined GW invariants from this (perturbed) moduli space.

\smallskip

We can follow as similar procedure for the family invariants by
introducing an inhomogeneous term into the $J_{\a}$-holomorphic
equation and vary $\nu$. This alternative definition of invariants
is more  geometric. In particular, using this definition of
invariants we can follow the analytic arguments of Ionel and
Parker  in \cite{ip2,ip3} to show sum formulas (\ref{E:Sum1}) and
(\ref{E:Sum2}) for the case at hand:  the class $s+df$ in $E(n)$.

\medskip

To simplify  notation in this section we will set $A=s+df$.

\bigskip
Using Prym structures defined as in \cite{lo},
we can lift the Deligne-Mumford space $\overline{\M}_{g,k}$
to a finite cover
  \begin{equation}\label{D:GC}
     p_{\mu} :  \overline{{\mathcal M}}_{g,k}^{\mu}
                \to
                \overline{{\mathcal M}}_{g,k}.
  \end{equation}
This finite cover is now a smooth manifold
and has a universal family
  \begin{equation*}
     \pi_{\mu} : \overline{{\mathcal U}}_{g,k}^{\mu}
                 \to
                 \overline{{\mathcal M}}_{g,k}^{\mu}
  \end{equation*}
which is projective.
Moreover, for each $b\in  \overline{{\mathcal M}}_{g,k}^{\mu}$,
$\pi_{\mu}^{-1}(b)$ is a stable curve isomorphic to
$p_{\mu}(b)$.

\smallskip
We fix, once and for all, an embedding of
$\overline{{\mathcal U}}_{g,k}^{\mu}$ into  some ${\Bbb{P}}^{N}$.
An inhomogeneous term $\nu$ is then defined as
a section of the bundle
$\mbox{Hom}(\pi^{*}_{1}(T{\Bbb{P}}^{N}),\pi^{*}_{2}TE(n))$
which is anti-$J$-linear :
  \begin{equation}
     \nu(j_{\Bbb{P}}(v))=-J(\nu(v))
     \ \ \ \ \ \ {\rm for\ any}\ v\in T{\Bbb{P}}^{N}
  \end{equation}
where $j_{\Bbb{P}}$ is the complex structure on ${\Bbb{P}}^{N}$.

\smallskip
For each stable map $f:C\to  E(n)$,
we can  specify one element $j\in p_{\mu}^{-1}(st(C))$.
Then $\pi_{\mu}^{-1}(j)$ is isomorphic to
the stable curve $st(C)$. In this way,
we can define a map
  \begin{equation}\label{D:phi}
     \phi : C \to  st(C) \,\cong \,
     \pi_{\mu}^{-1}(b) \,\subset\, \overline{{\cal U}}_{g,k}^{\mu}
     \hookrightarrow {\Bbb P}^{N}.
  \end{equation}
%

%%%%%  definition

\medskip
\begin{defn}
A stable $(J,\nu,\alpha)$-holomorphic map is a stable map
$f : (C,\phi) \to E(n)$
satisfying
  \begin{equation*}
     (\,df + J_{\alpha}dfj_{C}\,)(p)\ =\
     \nu_{\alpha}(\phi(p),f(p))
  \end{equation*}
where $\phi$ is defined as in (\ref{D:phi}),
and $\nu_{\alpha}=(I+JK_{\alpha})^{-1}\nu$.
\end{defn}

\bigskip
We denote the moduli space of stable $(J,\a,\nu)$-holomorphic maps
$(\, (f,(\phi,C),\alpha\,)$ by
  \begin{equation}\label{mstable}
     \CM_{g,k}(\,E(n),A,\nu,{\cal H},\mu\,)
  \end{equation}
where $\alpha$ in $\H$ and
$[f(C)]=A$ in $H_{2}(E(n);Z)$.
We also denote by
  \begin{equation*}
     \M_{g,k}(\,E(n),A,\nu,{\cal H},\mu\,)
  \end{equation*}
the set of
$\,( (f,(\phi,C)),\alpha \,)$
with a smooth domain $C$.
We will often abuse notation by writing $(f,C,\alpha)$, $(f,j,\alpha)$
or simply $(f,\alpha)$, instead of $(\,f,\,(\phi,C),\,\a\,)$.

\bigskip
There are stabilization and evaluation maps as in (\ref{st-ev}):
  \begin{equation}\label{E:st-ev}
     \M_{g,k}(\,E(n),A,\nu,{\cal H},\mu\,)\,
     @> st^{\mu}\times ev^{\mu} >> \,
     \CM_{g,k}^{\mu}\times E(n)^{k}.
  \end{equation}
Its Frontier is defined to be the set
  \begin{equation*}
     \{\,r\in\CM_{g,k}^{\mu}\times E(n)^{k}\,|\,
     r=\mbox{lim}(st^{\mu}\times ev^{\mu})
     (f_{n},\alpha_{n})\
     \mbox{and}\ (f_{n},\alpha_{n})\
     \mbox{has\ no\ convergent\ subsequence}\,\}.
  \end{equation*}

We denote by ${\cal Y}_{0}$ the space of all $\nu$
with $|\nu|_{\infty}$ is sufficiently small.
The following is the "Structure Theorem"
for the moduli space.

%%%%%   Structure Theorem %%%%%

\medskip
\begin{theorem}{\bf (\,Structure Theorem\,)}
                   \label{T:Structure}
For generic $\nu\in {\cal Y}_{0}$,  the space $\M_{g,k}(\,E(n),A,\nu,{\cal H},\mu\,)$
               is an smooth oriented manifold of dimension
                 \begin{equation}\label{E:dimMspace}
                   2\,c_{1}(A)\ +\  2\,(g-1)\  +\
                   2\,k\  +\  {\rm dim}(\H)\  =\  2\,(\,g\,+\,k\,).
                 \end{equation}
Furthermore, the frontier of the smooth map (\ref{E:st-ev})
               lies in dimension at most 2 less than $2(g+k)$.
\end{theorem}

\medskip
\noindent
{\bf Sketch of Proof.}\ \
The proof of this theorem is similar to
that of Proposition 2.3 in \cite{rt2}.
The  first statement follows from the standard argument
using Sard-Smale Theorem. To prove the second statement,
we first consider the well-defined stabilization and evaluation map
  \begin{equation}\label{E:st-ev-s}
     \CM_{g,k}(\,E(n),A,\nu,{\cal H},\mu\,)\,
     @> st^{\mu}\times ev^{\mu} >> \,
     \CM_{g,k}^{\mu}\times E(n)^{k}.
  \end{equation}
It then follows from Gromov Convergence Theorem \cite{is,p,pw} and
Lemma~\ref{L:UB-F} below that the stable moduli space
(\ref{mstable}) is compact and hence (\ref{E:st-ev}) extends
(\ref{E:st-ev-s}) continuously.

As in \cite{rt1,rt2}, we reduce the moduli space by
(i)  collapsing all ghost bubbles,
(ii) replacing each multiple maps from a bubble
     by its reduced map, and
(iii) identifying those bubble components which have the same image.
We denote this reduced moduli space by
  \begin{equation*}
     \CM_{g,k}^{r}(\,E(n),A,\nu,{\cal H},\mu\,).
  \end{equation*}
The map (\ref{E:st-ev}) now descends to the reduced moduli space
and by definition we have
  \begin{equation*}
     \mbox{Fr}(\,st^{\mu}\times ev^{\mu}\,)\,\subset\,
     st^{\mu}\times ev^{\mu}\left(\,
     \CM_{g,k}^{r}(\,E(n),A,\nu,{\cal H},\mu\,)\,
     \setminus\,\M_{g,k}(\,E(n),A,\nu,{\cal H},\mu\,) \,\right).
  \end{equation*}

It remains to show that  those strata consisting of $(f,\alpha)$
with domain more than two components has a dimension
at least 2 less than $2(g+k)$.
Similarly to the moduli space of $(J,\nu)$-holomorphic maps,
the strata corresponding to the domain with no bubble component has
a dimension at least 2 less than $2(g+k)$ for generic $\nu$.

On the other hand, it follows from
compactness of stable moduli space (\ref{mstable}) and
Theorem 2.4 of \cite{l} that
the restriction of $(f,\alpha)$ to any component of domain
should represents one of the following homology classes
  \begin{equation*}
     s,\ \ s+d_{1}f,\ \ d_{2}f \ \ \ \mbox{with}\ \ \
   0\ <\ d_{1},d_{2}\ \leq d
  \end{equation*}
Since the inhomogeneous term $\nu$ vanishes on bubble components,
by Theorem 2.4 of \cite{l} that each bubble component
maps into either a section or a singular fiber.

Now, suppose $(f,\alpha)$ has some bubble components.
Again by Theorem 2.4 of \cite{l}
either $\alpha\equiv 0$ or the zero divisor $Z(\alpha)$
contains some singular fibers.
Since there's no fixed component in the complete linear system
of a canonical divisor of $E(n)$, the parameter $\alpha$ lies in the proper
subspace of ${\cal H}$. This reduces the dimension of the strata
containing $(f,\alpha)$ at least 2.
\QED

%%%%%%%%%%%%  Definition of Invariants   %%%%%%%%%%%%

\vskip 0.8 cm Now, we are ready to define invariants. Instead of
using intersection theory as in \cite{rt1,rt2}, we will follow the
approach in \cite{ip2}. The above Structure Theorem and
Proposition 4.2 of [KM2] assert that the image
  \begin{equation*}
     st^{\mu}\times ev^{\mu}\left(\,
     \M_{g,k}(\,E(n),A,\nu,\H,\mu\,)\,\right)
  \end{equation*}
gives rise to a rational homology class in
$H_{*}(\, \CM_{g,k}^{\mu} ; {\Bbb Q} \,)\otimes
     H_{*}(\, E(n)^{k} ; {\Bbb Q}\, )$.
We denote it by
  \begin{equation}\label{D:homology}
     [\,\CM_{g,k}(\,E(n),A,\nu,\H,\mu\,)\,].
\end{equation}

%%%  Definition %%%

\medskip
\begin{defn}\label{def-fc}
For  $2g + k \geq 3$, we define invariants by
  \begin{equation*}
     GW_{g,k}(E(n),A,\H)(\beta;\a)\ =\
     \frac{1}{\lambda_{\mu}}\,
     (\,\beta\otimes \a\,)\ \cap\
     [\,\CM_{g,k}(\,E(n),A,\nu,\H,\mu\,)\,]
  \end{equation*}
where $\beta$ in $H^{*}(\,\CM_{g,k} ; {\Bbb Q}\,)$,
$\a$ in  $H^{*}(\,E(n)^{k} ; {\Bbb Q}\,)$, and
$\lambda_{\mu}$ is the order of
the finite cover in (\ref{D:GC}).
\end{defn}

%%%%% Comparison of Ruan-Tian and Li-Tian %%%%

\medskip

By repeating the same arguments for ordinary GW invariants,
we can show that these invariants are same as
the family invariants defined in Definition~\ref{D:fgw-inv't},
namely
  \begin{equation*}
    GW_{g,k}(E(n),A,\H) = GW_{g,k}^{\H}(E(n),A).
  \end{equation*}
In the below, we will not distinguish two invariants and
use the same notation
$GW_{g,k}^{\H}(E(n),A)$ for them.
We end this section by showing $F_{0}(0)=1$
which provides the initial condition for (\ref{ODE}).

\medskip
\begin{prop} \label{F(0)=1}
  \ \ $GW_{s,0}^{\H}(E(n))(f^{3})\ =\ 1$.
\end{prop}

\pf  Fix $\nu=0$. Since the section class $s$ is of type $(1,1)$,
Theorem 2.4 of \cite{l}
implies that  for any $(J,\alpha)$-holomorphic map
$(f,\alpha)$ with $[f]=s$,
$f$ is holomorphic and $\alpha=0$.
In fact, there is a unique such $f$ since $s^{2}=-n$.

Now, consider the linearization of $(f,\alpha)$-holomorphic
equation $L_{f}\oplus J\,df \oplus L_{0}$ as in appendix of
\cite{l}. Propositions A.1 and A.2 of the appendix of \cite{l}
show, quite generally, that $L_f$ is a $\dbar$ operator and
$L_{0}$ defines a map
  \begin{equation*}
     L_{0}:\H\to \mbox{Coker}(L_{f}\oplus J\,df)
  \end{equation*}
which is injective if and only if the family moduli space
$\CM^{\H}_{g,k}(E(n),A)$ is compact.  But we just showed the
moduli space is a
single point, and hence compact.

\smallskip
On the other hand,
$\mbox{Ker}(L_{f}\oplus Jdf)$ is same as
$H^{0}(f^{*}N)$, where
$N$ is the normal bundle of the section in $E(n)$.
It is trivial since the Chern number of $N$ is $s\cdot s\, =\, -n < 0$.
Therefore,
  \begin{equation*}
     \mbox{dim\,Coker}(\,L_{f}\oplus Jdf\,)
     \ =\  -\mbox{Index}(L_{f} \oplus Jdf )\ =\
     -2\,(\,c_{1}(f^{*}TE(n))- 1\,)\ =\  2\,(n-1)
  \end{equation*}
Since $L_{0}$ is injective and $\mbox{dim}(\H)=2(n-1)$,\
$L_{f}\oplus Jdf \oplus L_{0}$ is onto. That implies $\nu=0$
is generic
in the sense of Theorem~\ref{T:Structure}.
Consequently, the invariant is $\pm 1$.
In this case, the sign is determined by $L_{f}$ and
$L_{f}$ is $\dbar $-operator,
the invariant is 1.
\QED

\vskip 1cm

%%%%%%%%%%%%%%%%%%%%%%%%%%%%%%%%%%%%%%%%%%%%%%%%%%%%%%%%%%%%
%%%%%%%%%%%%%%%%%  Section 5 %%%%%%%%%%%%%%%%%%%%%%%%%%%%%%
%%%%%%%%%%%%%%%%%%%%%%%%%%%%%%%%%%%%%%%%%%%%%%%%%%%%%%%%%%%%%%

\setcounter{equation}{0}
\section{Degeneration of $E(n)$ }
\label{section5}
\bigskip

In this section, we describe a  degeneration
of $E(n)$ into a singular surface which is a union of
$E(n)$ and $E(0)$ with $V=T^{2}$ intersection.
We then define
the parameter space and inhomogeneous terms
corresponding to this degeneration.
The sum formulas (\ref{E:Sum1}) and (\ref{E:Sum2})
will be formulated from this degeneration

\bigskip
Let $D\subset {\Bbb C}$ be a small disk and
choose a smooth fibre $V$ in $E(2)$.
We denote by
  \begin{equation}\label{blow-up}
     p:Z\,\to \,E(n)\times D
  \end{equation}
the blow-up of $E(n)\times D$ along $V\times \{0\}$ and define\
%
  %\begin{equation}\label{degen}
     $\lambda: Z \,@> p >> \,E(n)\times D \,\to\, D$
  %\end{equation}
to be the composition map, where the second map is the projection
onto the second factor. The central fiber $Z_{0}=\lambda^{-1}(0)$
is a singular surface $E(n)\,\cup_{V}\,E(0)$ and the fiber
$Z_{\l}$ with $\l\ne 0$ is isomorphic to $E(n)$ as a complex
surface.

\smallskip
To save notation, we will use the same notation $(\omega,J,g)$ for
the induced K\"{a}hler structure on $Z$ and its restriction to
$Z_{\lambda}$, $E(n)$, and $E(0)$.

\smallskip
Fix a normal neighborhood $N_{E(n)}$ of $V$ in $E(n)$.
It is then a product $V\times D^{\prime}$,
where $D^{\prime}\subset {\Bbb C}$ is some disk.
Let $x$ be the holomorphic coordinate of $D^{\prime}$.
Then, the normal neighborhood $N$ of $V$ in $Z$ is given by
 \begin{equation*}%\label{lcoord}
  N\ =\ \{\ \left(\,v,x,\lambda,[\,l_{0};l_{1}]\,\right)\ |\
        v\in V,\ x\,l_{1} = \lambda\, l_{0}\ \}\
   \subset \ N_{E(n)}\times D^{\prime} \times D \times {\Bbb CP}^{1}
 \end{equation*}
where $[l_{0};l_{1}]$ is the homogeneous coordinates
of ${\Bbb CP}^{1}$.
It is covered by two patches
%\begin{equation*}
$U_{0}=(l_{0}\ne 0)$ and  %\mbox{\ \ \ and\ \ \ }
$U_{1}=(l_{1}\ne 0)$.
%\end{equation*}
On $U_{0}$, we set $y=l_{1}/\,l_{0}$.
Then we have
\begin{equation*}
N\ =\ \{\ (v,x,y)\ |\ v\in V\ \} \mbox{\ \ \ \ \ with\ \ \ \ \ }
\lambda(v,x,y) = xy.
\end{equation*}
Clearly, $Z_{\lambda}\,\cap\,N$ is given by
the equation $xy=\lambda$.
Note that we can also think of $y$ as a holomorphic normal coordinate
of the normal neighborhood $N_{E(0)}$ of $V$ in $E(0)$.

\begin{defn}
For some $\delta>0$ and $|\l|$,
we decompose $Z_{\l}$ as
a union of three pieces, two sides and a neck.
The {\em $\delta$-neck} is defined as
\begin{equation*}%\label{neck}
 Z_{\l}(\delta)\ =\ \left\{\ (v,x,y)\ \in\ Z_{\l}\,\cap\,N\ |\
                         \left|\,\,|x|^{2}-|y|^{2}\,\right|
                         \ \leq\ \delta\ \right\}.
\end{equation*}
$Z_{\l}\setminus Z_{\l}(\delta)$ consists of two components.
The $E(n)$-side is the component which contains
the region $|x|\, >\, |y|$,
while the component  $E(0)$-side contains the region $|x|\,<\,|y|$.
\end{defn}

On the neck region, there is a symplectic $S^{1}$-action with
Hamiltonian $t=\frac{1}{2}(\,|y|^{2}-|x|^{2})$. We can thus
decompose each $Z_{\l}$ as
 $Z_{\l} = Z_{\l}^{-}\,\cup\,Z_{\l}^{+}$,
where $Z_{\l}^{-}$ is a union of $E(n)$-side and the part of
$Z_{\l}(\delta)$ with $t\leq 0$. In fact, $E(n)\,(\mbox{resp}.\
E(0)\,)$ is the symplectic cut of $Z_{\l}^{-}\,(\mbox{resp}.\
Z_{\l}^{+})$ at $t=0$. Therefore, we have a collapsing map
\begin{equation}\label{collapsing}
 \pi_{\l}:Z_{\l}\to Z_{0}
\end{equation}
(cf. section 2 of \cite{ip3}).

\vskip 0.6cm

Next, we define the parameter spaces. Let $U$ be a neighborhood of
$V$ in $E(n)$ that does not contain any singular fibers. Choose a
bump function $\beta$ which satisfies $\beta = 1$ on
$E(n)\setminus U$ and $\beta=0$ near $V$ in $U$.

\begin{defn}\label{para}
We define the parameter spaces by
  \begin{equation*}
     \H_{\l}\ =\ \{\ \,\a_{\l} =\, p^{*}_{\l}\,\beta\,\a\ \,|\,\
                        \a\in \H\ \,\}
     \ \ \ \ \
     \mbox{and}\ \ \ \ \
     \H_{E(n)}\ =\ \{\ \beta\,\a\ \,|\,\
                        \a\in \H\ \,\}
  \end{equation*}
where $p_{\l}$ is the restriction of (\ref{blow-up}) to $Z_{\l}$.
\end{defn}

Note that $\H_{E(n)}=\{0\}$ for $n=0,1$. On the other hand, each
$\a\in \H_{E(n)}$ (\,resp. $\a_{\l}\in \H_{\l}$\,) is
$J$-anti-invariant and $\a=0\,(\mbox{resp.}\ \a_{\l}=0)$ near $V$
by definition. Hence $J_{\a}=J$\,(resp. $J_{\a_{\l}}=J$) near $V$.

\bigskip
Lastly, following \cite{ip3}, we define inhomogeneous terms. An
inhomogeneous term $\nu$ of the the fibration $\l:Z\to D$ is a
section of the bundle $\mbox{Hom}(T{\Bbb P}^{N},TZ)$ over ${\Bbb
P}^{N}\times Z$ for some ${\Bbb P}^{N}$, which satisfies
Definition 2.2 of \cite{ip3}. We denote by $\J_{0}(Z)$ the space
of all such $\nu$ with sufficiently small $|\nu|_{\infty}$ and use
the same notation $\nu$ for the restriction of $\nu$ to $Z_{\l}$,\
$E(n)$,\ and $E(0)$.

%Lastly, following \cite{ip2}, we define inhomogeneous terms. Let
%$N_{V}$ be the normal bundle of $V$ in $Z$ and $\xi\to \xi^{N}$ be
%the orthogonal projection onto $N_{V}$.

%\begin{defn}
%After noting that $J$ is K\"{a}hler, we define an inhomogeneous
%term $\nu$ of $Z$ as a section of $\mbox{Hom}(T{\Bbb P}^{N},TZ)$
%over ${\Bbb P}^{N}\times Z$ for some ${\Bbb P}^{N}$, which
%satisfies
%
%  \begin{equation*}
%    \nu^{N} = 0    \mbox{\ \ \ and\ \ \ }
%    [\,\nabla_{\xi}\nu + J\nabla_{J\xi}\nu\,]^{N}=0
%  \end{equation*}
%for all $\xi\in N_{V}$ (cf. Definition 3.2 of \cite{ip2} and
%Definition 2.2 of \cite{ip3}).
%\end{defn}

%We denote by
%$\J_{0}(Z)$
%the space of all such $\nu$ with sufficiently small $|\nu|_{\infty}$
%and use the same notation $\nu$ for the restriction
%of $\nu$ to $Z_{\l}$,\ $E(n)$,\ and $E(0)$.

\vskip 1cm

%%%%%%%%%%%%%%%%%%%%%%%%%%%%%%%%%%%%%%%%%%%%%%%%%%%%%%%%%%%%%%%%%%%%%%%%%%%
%%%%%%%%%%%%%%%%%%%%%%%%%%%%%%%%%%%%%%%%%%%%%%%%%%%%%%%%%%%%%%%%%%%%%%%%%%%
\setcounter{equation}{0}
\section{\bf Splitting of Maps }
\label{section6}
\bigskip

In this section, we show the uniform energy bound of maps and
$L^{2}$-bound of parameters $\a$. By Gromov Convergence Theorem,
these leads to the compactness of family moduli spaces. The
splitting arguments as in section 3 of \cite{ip3} then follows
from the compactness and the choice of inhomogeneous terms and
parameters $\a$ --- we bumped $\a$ to 0 along $V$.

\medskip
Let $(X,\omega,h,J)$ be a 4-dimensional almost K\"{a}hler
manifold. Recall that $\a$ is a $J$-anti-invariant 2-form on $X$
if $\a(Ju,Jv)=-\a(u,v)$. Fix a metric within the conformal class
$j$ on a Riemann surface $(C,j)$ and let $dv$ be the associated
volume form.

\begin{lemma}\label{Cor1.4}
For any $C^{1}$ map $f:C\to X$  if $\a$ is $J$-anti-invariants,
then
\begin{equation}\label{easy-forget}
  f^{*}\a\ \leq\ 2\,|\a|\,|df|\,|\dbar_{J} f|.
\end{equation}
On the other hand, if the map $f$ is $(J,\nu,\a)$-holomorphic,
then we have
\begin{align}
     |\dbar_{J} f|^{2}\,dv\  & =\  f^{*}\a  +
     2\,\langle \dbar_{J} f , \nu \rangle\,dv
     \label{E:EnergyB} \\
     (\,1 +f^{*}(|\a|^{2})\,)\,f^{*}\omega\ & =\
     \frac{1}{2}\,(\,1-f^{*}(|\a|^{2})\,)\,|df|^{2}\,dv  -
     4\,\langle \dbar_{J} f , \nu \rangle\,dv +  4|\nu|^{2}\,dv.
     \label{E:alphaB}
\end{align}
\end{lemma}

\pf The proof of (\ref{E:EnergyB}) and (\ref{E:alphaB}) is similar
to those of Corollary 1.4 of \cite{l}. We will prove
(\ref{easy-forget}) only. Fix a point $z\in C$ and an orthogonal
basis $\{e_{1},e_{2}=je_{1}\}$ of $T_{z}C$. Then we have
\begin{align}\label{easy-forget1}
  \a\big(\,df(e_{1}),df(e_{2})\,\big)\ &=\
  \a\big(\,df(e_{1}),df(e_{2})+Jdf(je_{2})-Jdf(je_{2})\,\big)\notag \\
  &=\
  \a\big(\,df(e_{1}),2\dbar_{J}f(e_{2})\,\big)\ +\
  \a\big(\,df(e_{1}),Jdf(e_{1})\,\big).
\end{align}
Since $\a$ is $J$-anti-invariant,
$\a\big(\,df(e_{1}),Jdf(e_{1})\,\big)=0$. Therefore,
(\ref{easy-forget}) follows from (\ref{easy-forget1}). \QED

\bigskip
We denote the stable family moduli space of
$(J,\nu,\a_{\l})$-holomorphic maps $(f,\a_{\l})$ by
  \begin{equation*}
     \CM_{g,k}(\,Z_{\l},s+df,\nu,\H_{\l}\,),
     \mbox{\ \ \ or\ \ simply\ \ \ }
     \CM_{g,k}(\,\l,d\,)
  \end{equation*}
where $\nu\in \J_{0}(Z)$  and $\a_{\l}\in \H_{\l}$. Compactness of
the family moduli space follows from Gromov Convergence Theorem
and the following lemma.

\begin{lemma}\label{L:UB-F}%%Uniform Bound for Fibration
Let $|\nu|_{\infty}$ be sufficiently small. Then, there exit
uniform constants $E_{d}$ and $N$, which does not depend on $\l$,
such that
  \begin{equation*}
     E(f)\ =\ \frac{1}{2}\,\int_{C}|df|^{2}\ \leq\ E_{d}
     \mbox{\ \ \ \ and\ \ \ \ }
     ||\a_{\l}||\ =\ \int_{Z_{\l}}\a_{\l}\wedge\a_{\l}\  \leq\ N
  \end{equation*}
for any $(f,C,\a_{\l})$\  in\ $\CM_{g,k}(\,\l,d\,)$.
\end{lemma}

\pf
We first show uniform bound of $||\a_{\l}||$. This proof is
similar to those of Lemma 4.4 except for using (\ref{E:alphaB})
instead of Corollary 1.4b of \cite{l}. For each $\a_{\l}$ in
$\H_{\l}$, we choose a sufficiently small neighborhood of
$N(\a_{\l})$ of the zero set of $\a_{\l}$ and let $m(J_{\l})$ and
$N$ be as in the proof of Lemma 4.4 of \cite{l}. If there is  a
holomorphic fiber $F\subset E(n)\setminus N(\a_{\l})$ such that
  \begin{enumerate}
     \item[(i)]
        $f$ is transversal to  $F$,
     \item[(ii)]
        at each $p\in f^{-1}(F)$, $f$ is transversal to a holomorphic
        disk $D_{f(p)}$ normal to $F$
        at $f(p)$, and
     \item[(iii)]
        $4\,|df|\,|\nu|\ +\ 4\,|\nu|^{2}\ \leq\
        \frac{1}{2}\,|df|^{2}$ on $f^{-1}(F)$
  \end{enumerate}
then the proof follows exactly as in the proof of Lemma
Lemma 4.4 of \cite{l}.  We can clearly find
fibers satisfying (i) and (ii),
so  we need only verify that we can
also obtain (iii).
For that we consider the set $C_{0}$ of all points in
$C$ where\ \,
$4\,|df|\,|\nu| +  4|\nu|^{2} \, >\,
\frac{1}{2}\,|df|^{2}$.
Then $|df|^{2}\leq 100|\nu|^{2}$ on $C_0$.
Therefore
  \begin{equation}\label{E:Area}
     \int_{C_{0}}|\,d\pi\circ df\,|^{2} \ \leq \
     100\,\mbox{Area}(\,st(C)\,)\,
     |\,d\pi\,|_{\infty}^{2}|\,\nu\,|^{2}_{\infty}
  \end{equation}
where $\pi :Z_{\l}\to {\Bbb CP}^{1}$ is the elliptic structure for
$J$ on $Z_{\l}$. We can thus assume that the left hand side of
(\ref{E:Area}) is less than $\frac{1}{3}\mbox{Area}({\Bbb
CP}^{1})$ for sufficiently small $|\nu|_{\infty}$. On the other
hand, from the definition of $N(\a_{\l})$, we can also assume that
$\mbox{Area}(\,\pi( N(\a_{\l}))\, )\leq
\frac{1}{3}\mbox{Area}({\Bbb CP}^{1})$. Therefore, we can always
choose a holomorphic fiber $F = \pi^{-1}(q)$ as in the above claim
with $q\in {\Bbb CP}^{1}\setminus (\  \pi(N(\a_{\l}))\cup \pi\circ
f(C_{0})\,)$.

\smallskip Next, we show uniform bound of the energy $E(f)$. By
definition~\ref{para},  $\a_{\l}=p^{*}_{\l}\,(\beta\,\a)$ for some
$\a\in\H$ and $p^{*}_{\l}(\a)$ is $J$-anti-invariant. We define
$C_{-}$ as the set of all $z$ in $C$ with
$f^{*}\,p^{*}_{\l}\,\a(\,e_{1}(z),e_{2}(z)\,)\leq 0$, where
$\{\,e_{1}(z),e_{2}(z)=j\,e_{1}(z)\,\}$ is an orthonormal basis of
$T_{z}C$. Then (\ref{E:EnergyB}) implies that $ |\dbar_{J}f|\leq
|\nu|$ on $C_{-}$ and hence by (\ref{easy-forget}) we have
  \begin{equation}\label{mid}
     0\ \leq\
      -f^{*}\,p^{*}_{\l}\,\a(\,e_{1}(z),e_{2}(z)\,)\
      \leq\ 2|p^{*}_{\l}\,\a|\,|df|\,|\dbar_{J} f|\
      \leq\ 2\,M\,|df|\,|\nu|\
  \end{equation}
for any $z\in C_{-}$, where $M\,=\,\mbox{max}\{\
|p_{\l}^{*}\alpha|\ |\ ||\a_{\l}||\,\leq N\ \}$. Therefore, we can
conclude that
\begin{align*}
\frac{1}{2}\int_{C}|df|^{2}\  & =\ \int_{C}|\dbar_{J}f|^{2} \ +\
\omega(s+df)\\
& \leq
    \int_{C}f^{*}\,p^{*}_{\l}(\beta\alpha) +
    2\int_{C}|df||\nu| + \omega(s+df)\\
& \leq
    \int_{C\setminus C_{-}}f^{*}p^{*}_{\l}\alpha +
    2\int_{C}|df||\nu| + \omega(s+df) \\
& \leq
    -\int_{C_{-}}f^{*}p^{*}_{\l}\alpha +
    2\int_{C}|df||\nu| + \omega(s+df) \\
& \leq
    (1 + 2M)\left( \int_{C}|\nu|^{2} \right)^{\frac{1}{2}}
           \left( \int_{C}|df|^{2} \right)^{\frac{1}{2}}
     + \omega(s+df)
\end{align*}
where the second inequality follows from (\ref{E:EnergyB}), the
fourth inequality follows from  $p_{\l}^{*}\a\,(s+df)=0$ and the
last from (\ref{mid}). This implies the uniform energy bound
independent of $\lambda$ for sufficiently small $|\l|_{\infty}$.
\QED

\begin{rem}\label{R:AI=RI}
Repeating the same argument as in section 4, one can show that the
moduli spaces
  \begin{equation*}
    \CM_{g,k}(E(n),s+df,\H_{E(n)},\nu)
    \mbox{\ \ \ and\ \ \ }
    \CM_{g,k}(Z_{\l},s+df,\H_{\l},\nu)
  \end{equation*}
defines family invariants $GW_{g,k}(E(n),A,\H_{E(n)})$ and
$GW_{g,k}(Z_{\l},A,\H_{\l})$, respectively.
Moreover, by the standard corbodism argument as in Lemma 4.9 of \cite{rt2}
\ we have
  \begin{equation}\label{cobo}
     GW_{g,k}(E(n),s+df,\H_{E(n)})=GW^{\H}_{g,k}(E(n),s+df)=
     GW_{g,k}(Z_{\l},s+df,\H_{\l}).
  \end{equation}
\end{rem}

\bigskip
The following shows how maps into $Z_{\l}=E(n)$ split along the
degeneration of $E(n)$. It is also a key observation for gluing of
maps into $E(n)$ and $E(0)$, which leads to the sum formulas
(\ref{E:Sum1}) and (\ref{E:Sum2}).

\begin{lemma}\label{L:splitting}
Let $\nu\in\J_{0}(Z)$ and $\{(f_{\l},C_{\l},\alpha_{\l})\}$ be any
sequence with $(f_{\l},C_{\l},\alpha_{\l})\in
\CM_{g,k}(\,\l,d\,)$. Then as $\l\to 0$,\ $f_{\l}$ converges to a
limit
           $f_{0}:C_{0}\to Z_{0}$\ %=E(n)\cup_{V} E(0)$
           and $\a_{\l}$ converges to $\a_{0}$,
           after passing to some subsequences, such that
\begin{enumerate}
\item[(a)] the limit map $f_{0}$ can be decomposed as
              \begin{equation*}
                f_{1}:C_{1}\to E(n),\ \
                f_{2}:C_{2}\to E(0),\ \
                \mbox{and\ \ }
                f_{3}:C_{3}\to V
               \end{equation*}
             and
             $f_{1}$\,(\,resp. $f_{2}$\,) represents homology class
             $s+d_{1}f$ in $E(n)$\,(\,resp. $s+d_{2}f$ in E(0)\,) and
             $f_{3}$ represents $d_{3}[V]$ in $V$ with
             $d_{1}+d_{2}+d_{3}=d$,
\item[(b)] for $i=1,2$, each $f_{i}$ transverse to $V$
              with $f_{i}^{-1}(V)=\{\,p_{i}\,\}$, where $p_{i}$
              is a node of $C$.
              %and,
%\item[(c)] for generic $\nu\in \J_{0}(Z)$ if no marked point of $C_{0}$
           %maps into $V$, then
           %the (arithemetic) genus of $C_{3}$ is at least 2.
\end{enumerate}
\end{lemma}

\pf
By Gromov Convergence Theorem and Lemma~\ref{L:UB-F},
$f_{\l}$ converges to a limit $f_{0}:C_{0}\to Z_{0}$.
Since $\a_{\l}=0$ near $V\subset Z$,\
we have $J_{\a_{\l}}=J$
near $V$ in $Z$.
Therefore, (a) and (b) follows from Lemma 3.4 of \cite{ip2} and
Lemma 3.3 of \cite{ip3}.
%jonghyun

\vskip 1cm

%%%%%%%%%%%%%%%%%%%%%%%%%%%%%%%%%%%%%%%%%%%%%%%%%%%%%%%%%%%%
%%%%%%%%%%%%%%%%%  Section 7 %%%%%%%%%%%%%%%%%%%%%%%%%%%%%%
%%%%%%%%%%%%%%%%%%%%%%%%%%%%%%%%%%%%%%%%%%%%%%%%%%%%%%%%%%%%%%

\setcounter{equation}{0}
\section{Relative Invariants of $E(n)$  }
\label{section7}
\bigskip

In this section, following \cite{ip2}, we define relative invariants
of $E(n)$ relative to a smooth elliptic fiber $V=T^{2}$.
In our case, the {\em rim tori} in $E(n)\setminus V$ disappear
when we glue $E(n)$ and $E(0)$ along $V$.
Together with the simple matching condition as in Lemma~\ref{L:splitting},
that observation leads to the simple definition of relative invariants.

\bigskip
As in section \ref{section4}, we fix the complex structure on
$E(n)$. We also  assume that we always work with a finite good
cover $p_{\mu}$ as in (\ref{D:GC}) without specifying it.
Throughout this section, $A$ always denotes the class $s+df$.

\smallskip
For $\nu$ in $\J_{0}(Z)$,
we define the relative moduli space by
\begin{equation*}\label{relM-X}
\M^{V}_{g,k+1}(E(n),A,\H_{E(n)},\nu)  =
\left\{\left(\,f,\alpha\,\right)
        \in \CM_{g,k+1}(E(n),A,\H_{E(n)},\nu)\ |\
        f^{-1}(V)=\{x_{k+1}\}\,
        \right\}.
\end{equation*}
As in \cite{ip2}, we compactify this moduli space
by taking its closure
  \begin{equation*}
     C\M^{V}_{g,k+1}(E(n),A,\H_{E(n)},\nu)\
  \end{equation*}
in the space of stable maps $\CM_{g,k+1}(E(n),A,\H_{E(n)},\nu)$.
Note that for each $\a\in \H_{E(n)}$,\ $\a=0$ in some neighborhood
of $V\subset E(n)$ and hence $J_{\a}=J$ on that neighborhood.
Therefore, Proposition~\ref{P:Structure-R} below follows from the
same arguments as in Lemma 4.2 and Proposition 6.1 of \cite{ip2},
and Theorem~\ref{T:Structure}.

\begin{prop}\label{P:Structure-R}
   For generic $\nu$ in $\J_{0}(Z)$
\begin{enumerate}
\item[(a)] $\M_{g,k+1}^{V}(\,E(n),A,\H_{E(n)},\nu\,)$
             is an orbifold of dimension $2+2(g+k)$
             for $n=0$ and $2(g+k)$ for $n\geq 1$, and
\item[(b)] the Frontier of the map
    \begin{equation}\label{E:F-RIVT}
             {\cal M}_{g,k+1}^{V}(E(n),A,\H_{E(n)},\nu)\,
             @>\  st\times ev\times h\  >> \,
             \overline{{\cal M}}_{g,k+1}\times E(n)^{k}\times V
    \end{equation}
            is contained in codimension at least 2,
            where $ev$ is the evaluation map of the first $k$ marked
            points and $h$ is the evaluation map of the last marked point.
\end{enumerate}
\end{prop}

\medskip
Proposition~\ref{P:Structure-R} together with
Proposition 4.2 of \cite{km2} assert that
the image of
(\ref{E:F-RIVT}) gives rise to a rational homology class.
We denote it by
\begin{equation*}
[\,\M_{g,k+1}^{V}(E(n),A,\H_{E(n)})\,]\in
H_{*}(\overline\M_{g,k+1};{\Bbb Q})\otimes
H_{*}(E(n)^{k};{\Bbb Q})\otimes H_{*}(V;{\Bbb Q}).
\end{equation*}

\smallskip
\begin{defn}\label{D:RI}
For $2g+k\geq 3$, we define relative invariants  by
\begin{equation*}
GW^{V}_{g,k+1}(E(n),A) (\,\beta;\alpha;C(\gamma)\,)=
\left(\,\beta\otimes\alpha\otimes \gamma\,\right) \cap
[\,\M_{g,k+1}^{V}(E(n),A,\H_{E(n)})\,]
\end{equation*}
where $\beta \in H^{*}(\overline\M_{g,k+1};{\Bbb Q})$,\
$\alpha \in H^{*}(E(n)^{k};{\Bbb Q})$,\  and
$\gamma \in H^{*}(V;{\Bbb Q})$.
\end{defn}

\medskip
The relative invariants of $E(0)$ and $E(1)$ defined as in
Definition~\ref{D:RI} are less finer than those in \cite{ip2}
(\,cf. Appendix in \cite{ip3}\,).  On the other hand, for another
smooth fiber $U=T^{2}$ of $E(0)$ we can define relative invariants
relative to both $V$ and $U$ as in Definition~\ref{D:RI}.  In the
below, we will denote ordinary and relative GW invariants for
$E(0)$ by
\begin{equation*}
\Phi_{A,g}\, ,\ \ \  \Phi^{V}_{A,g}\, ,\ \ \ \mbox{and}\ \ \
\Phi^{V,U}_{A,g}\, ,\ \ \mbox{respectively}.
\end{equation*}

\bigskip
We end this section by relative invariants of $E(0)$ for the class
$s+df$. Recall that for positive integer $d$,\ $\sigma(d)$ is the
sum of the divisors, namely $\sigma(d)=\sum_{k|d}\,k$. For
convenience we set $\sigma(0)=-1/24$.

\begin{lemma}{\rm (\cite{ip3})}\label{rel-compu}\
Let $V\subset E(0)$ be a smooth elliptic fiber.
  \begin{enumerate}
     \item[(a)]
        $\Phi^{V}_{s+df,0}(\,\tau(f^{*});C(f)\,) = 0$.
     \item[(b)]
        $\Phi^{V}_{s+df,1}(\,\tau(f^{*});C(pt)\,) = 2\sigma(d)$.
     \item[(c)]
        $\Phi^{V}_{s+df,0}(\,pt;C(f)\,)=
        \Phi^{V}_{s+df,0}(\,C(pt)\,)=
        1$ if $d=0$ and 0 otherwise.
     \item[(d)]
        $\Phi^{V}_{s+df,1}(\,pt;C(pt)\,)= d\,\sigma(d).$
     \item[(e)]
        $\Phi^{V,U}_{s+df,1}(\,C(pt),C(pt)\,)=0$.
  \end{enumerate}
\end{lemma}

\vskip 1cm

%%%%%%%%%%%%%%%%%%%%%%%%%%%%%%%%%%%%%%%%%%%%%%%%%%%%%%%%%%%%
%%%%%%%%%%%%%%%%%  Section 8 %%%%%%%%%%%%%%%%%%%%%%%%%%%%%%
%%%%%%%%%%%%%%%%%%%%%%%%%%%%%%%%%%%%%%%%%%%%%%%%%%%%%%%%%%%%%%

\setcounter{equation}{0}
\section{Sum Formula  }
\label{section8}
\bigskip

This section shows the sum formulas (\ref{E:Sum1}) and
(\ref{E:Sum2}) using a family version of Gluing Theorem --- a map
into $Z_{0}$ satisfying the matching condition as in
Lemma~\ref{L:splitting} can be smoothed to produce a map into
$Z_{\l}$. That smoothing relates invariants of $Z_{\lambda}=E(n)$
with relative invariants of $E(n)$ and $E(0)$ relative to a smooth
fiber $V=T^{2}$.

\medskip
Throughout this section we fix $n\geq 1$.  Recall the evaluation
map of last marked point as in (\ref{E:F-RIVT}). There is an
evaluation map
  \begin{equation*}
    ev_{V}:
    \bigcup\,(\  \M^{V}_{g_{1},k_{1}+1}(\,E(n),s+d_{1}f,\H_{E(n)},\nu)
     \times\,
     \M^{V}_{g_{2},k_{2}+1}(\,E(0),s+d_{2}f,\nu)\ )
    \to
     V^{2}
  \end{equation*}
which records the intersection points with $V$, where the union is
over all $g_{1}+g_{2}=g$, $k_{1}+k_{2}=k$ and $d_{1}+d_{2}=d$. We
set
  \begin{equation}\label{rel}
     \M^{V}_{g,k}(\,d\,) \ =\ ev_{V}^{-1}(\triangle)
  \end{equation}
where $\triangle$ is the diagonal of $V^{2}$. This space is an
orbifold of dimension $2(g+k)$ for generic $\nu$ in $\J_{0}(Z)$
and comes with stabilization and evaluation maps
\begin{equation*}
  \M^{V}_{g,k}(\,d\,)\ @>\  st\times ev\  >>\  \dgcc\ =\ \dgc
\end{equation*}
where the union is over all $g=g_{1}+g_{2}$, and $k=k_{1}+k_{2}$.

\bigskip
Now, consider a sequence of  maps $(f_{\l},\a_{\l})$ in
$\CM_{g,k}(\,\l,d\,)$. By Lemma~\ref{L:splitting}, as $\l\to 0$,
the maps $(f_{\l},\a_{\l})$ converge to a limit $(f_{0},\a_{0})$,
after passing to some subsequences. In general, the limit
$(f_{0},\a_{0})$ might be not in $\M^{V}_{g,k}(\,d\,)$. That
happens if some components of $f_{0}$ map entirely into $V$. On
the other hand, by Lemma 1.5 of \cite{ip2} there is a constant
$c_{V}$, depending only on $(J_{V},\nu_{V})$ such that every
stable $(J_{V},\nu_{V})$-holomorphic maps have energy great than
$c_{V}$. This implies that for small $|\l|$  the energy of
$f_{\l}$ in the $\delta$-neck
\begin{equation}\label{energy}
 E^{\delta}(f_{\l})\ =\ \frac{1}{2}\int |df_{\l}|^{2} + |d\phi|^{2}
\end{equation}
is greater than  $c_{V}$, where the integral is over
$f_{\l}^{-1}(\,Z_{\l}(\delta)\,)$ and $\phi:C_{\l}\to {\Bbb
P}^{N}$ as in (\ref{D:phi}). Therefore, for each $\l$ if $f_{\l}$
is $\delta$-flat\,(see Definition~\ref{flat} below) then $f_{0}$
is also  $\delta$-flat and hence the limit $(f_{0},\a_{0})$ is
contained in $\M^{V}_{g,k}(\,d\,)$.

\medskip
Following \cite{ip3}, we define  $\delta$-flat maps as follows:

\begin{defn}\label{flat}
A stable $(J,\nu,\alpha)$-holomorphic map $(f,\a)$ into $Z_{\l}$
is $\delta$-flat if
  \begin{equation}\label{flat-e}
     E^{\delta}(f)\ \leq\  \frac{c_{V}}{2}
  \end{equation}
\end{defn}

Note that any $\delta$-flat map $(f,\a)$ into $Z_{0}$ has no
component maps  into $V$. We denote by
  \begin{equation}
     \CM^{\,\delta}_{g,k}(\,\l,d\,)\ \subset\
     \CM_{g,k}(\,\l,d\,)
     \mbox{\ \ \ \ \ (\ resp.\ \ }
     \M^{V,\delta}_{g,k}(\,d\,)\ \subset\
     \M^{V}_{g,k}(\,d\,)\ )
  \end{equation}
the set of all $\delta$-flat maps in $\CM_{g,k}(\,\l,d\,)$
(\,resp. in $\M^{V}_{g,k}(\,d\,)$\,).

\bigskip The following is a family version of Theorem 10.1 of
\cite{ip3}. It shows that a $\delta$-flat map into $Z_{0}$ can be
smoothed to produce a $\delta$-flat map into $Z_{\l}$ for small
$|\l|$.

\begin{theorem}\label{T:Gluing}
For generic $\nu\in\J_{0}(Z)$ and for small $|\l|$,
there is a diagram
$$
\xymatrix{
 \dga \ar[rr]^{\Phi_{\l}} \ar[d]^{st\times ev}
   && \dgb   \ar[d]_{st\times ev} \\
 \dgcc \ar[dr]^{\sigma\times \pi_{0}}  && \dgd \ar[dl]_{id\times \pi_{\l}} \\
           & \dge      }
$$
which commutes up to homotopy, where $\Phi_{\l}$ is an embedding,
$\sigma$ is the gluing map of the domain as in (\ref{gluing-1}),
$\pi_{\l}$ is the collapsing map as in (\ref{collapsing}), and
$\pi_{0}\ :\ \bigcup\ (\, E(n)^{k_{1}}\times E(0)^{k_{2}}\,)\ \to\
                        Z_{0}^{k}$ defined by
%\begin{equation}\label{collap1}
%  \pi_{0}:\bigcup\ \big( E(n)^{k_{1}}\times E(0)^{k_{2}}\,big)\to
%                        Z_{0}^{k}
%\end{equation}
$\pi_{0}(x_{1},\cdots,x_{k_{1}},y_{1},\cdots,y_{k_{2}})=
        (x_{1},\cdots,x_{k_{1}},y_{1},\cdots,y_{k_{2}})$.
%where the union is over all $k=k_{1}+k_{2}$.
\end{theorem}

\bigskip
First, we use Theorem~\ref{T:Gluing} to derive the sum formula
(\ref{E:sumf}) for certain constraints. Let
$\b=\b_{1}\otimes\cdots\otimes\b_{k}\in H^{2r}(Z_{0}^{k})$, where
$r=g+k$. Denote by $B^{i}$ a geometric representative of the
Poincar\'{e} dual of $\b_{i}$. We assume that for some $0\leq
k_{1}\leq k$
\begin{enumerate}
\item[(i)] each $B^{i}$ lies in $E(n)$-side if $i\leq k_{1}$ and in
$E(0)$-side if $i>k_{1}$, and
\item[(ii)]
$\mbox{deg}(\b_{1}\otimes\cdots \otimes\b_{k_{1}})=2(g_{1}+k_{1})$
for some $0\leq g_{1}\leq g$.
\end{enumerate}
Note that the assumption implies that there is a decomposition
\begin{equation*}\label{constraint}
\pi_{0}^{*}\b\ =\ \b_{1}+\b_{2}\ \ \ \mbox{with}\ \ \ \b_{1}\in
H^{2(g_{1}+k_{1})}(\,E(n)^{k_{1}})\ \ \mbox{and}\ \ \b_{2}\in
H^{2(g_{2}+k_{2})}(\,E(0)^{k_{2}})
\end{equation*}
where $g_{2}=g-g_{1}$ and $k_{2}=k-k_{1}$. On the other hand, the
inverse image of $B^{i}$ under $\pi_{\l}$  gives a continuous
family of geometric representatives $B^{i}_{\l}$ of the
Poincar\'{e} dual of $\pi_{\l}^{*}\b_{i}$ in $H^{*}(Z_{\l})$. We
define the cut-down moduli spaces by
\begin{align}
\CM_{g,k}(\,\l,d\,)\,\cap\,\pi_{\l}^{*}\b\ &=\ \big\{\
(f,\a)\in\CM_{g,k}(\l,d)\ |\ ev_{i}(f,\a)\in B_{\l}^{i}\ \big\}
\label{cut-down1}\\
\M_{g,k}^{V}(\,d\,)\,\cap\,\pi_{0}^{*}\b \ &=\ \big\{\
\big(\,(f_{1},\a),f_{2}\big)\in\M_{g,k}^{V}(d)\ |\
ev_{i}\big(\,(f_{1},\a),f_{2}\big)\in B^{i}\
\big\}\label{cut-down2}
\end{align}
where $ev_{i}$ is the evaluation map of the i-th marked point.
Both cut-down moduli spaces (\ref{cut-down1}) and
(\ref{cut-down2}) are finite. In particular, any maps in
(\ref{cut-down1}) is $\delta$-flat for some $\delta>0$.

\medskip
\begin{prop}\label{P:sumf}
Let $\b\in H^{2r}(Z_{0}^{k})$ be a constraint as above. Then
\begin{align}\label{E:sumf}
GW^{\H}_{s+df,g}(\,\b\,)\ &= \sum_{d_{1}+d_{2}=d}
     GW^{V}_{s+d_{1}f,g_{1}+1}
     (\,\b_{1};C(pt)\,)\,
     \Phi^{V}_{s+d_{2},g_{2}-1}(\,\b_{2};C(f)\,)
     \notag \\
     &+
     \sum_{d_{1}+d_{2}=d}
     GW^{V}_{s+d_{1}f,g_{1}}(\,\b_{1};C(f)\,)\,
     \Phi^{V}_{s+d_{2}f,g_{2}}(\,\b_{2};C(pt)\,).
\end{align}
\end{prop}

\pf Denote the set of limits of sequences of maps in
(\ref{cut-down1}) as $\l\to 0$ by
\begin{equation}\label{limitset}
\underset{\l\to
0}{\lim}\left(\,\CM_{g,k}(\,\l,d\,)\,\cap\,\pi_{\l}^{*}\b\,\right).
\end{equation}
We first assume that the limit set (\ref{limitset}) is contained
in the space (\ref{rel}). Then for small $|\l|$ all maps in
(\ref{cut-down1}) should be $\delta$-flat. In that case
Theorem~\ref{T:Gluing} implies that
\begin{equation}\label{sumfm}
(id\times\pi_{\l})_{*}\big
[\,\CM_{g,k}(\,\l,d\,)\,\cap\,\pi_{\l}^{*}\b\,\big]\ =\
(\sigma\times\pi_{0})_{*}
\big[\,\M_{g,k}^{V}(\,d\,)\,\cap\,\pi_{0}^{*}\b\,\big]
\end{equation}
as a homology class in $H_{0}\big(\,\CM_{g,k}\times
Z_{0}^{k}\,;\,{\Bbb Q}\,\big)$. The left-hand side of
(\ref{sumfm}) becomes
\begin{equation}\label{right}
\big [\,\CM_{g,k}(\,\l,d\,)\,\cap\,\pi_{\l}^{*}\b\,\big]\,=\,
GW_{s+df,g}\big(Z_{\l},H_{\l}\big)\big(\b\big)\,=\,
GW^{\H}_{s+df,g}\big(E(n)\big)\big(\b\big)
\end{equation}
where the second equality follows from (\ref{cobo}). On the other
hand, by assumption on the constraint $\b$ and the routine
dimension count, the right-hand side of (\ref{sumfm}) becomes
\begin{align}\label{left}
     \big[\,\M_{g,k}^{V}(\,d\,)\,\cap\,\pi_{0}^{*}\b\,\big]\
     &=
     \sum_{d_{1}+d_{2}=d}
     GW^{V}_{s+d_{1}f,g_{1}+1}(\,\b_{1};C(pt)\,)\,
     \Phi^{V}_{s+d_{2},g_{2}-1}(\,\b_{2};C(f)\,)
     \notag \\
     \ &+
     \sum_{d_{1}+d_{2}=d}
     GW^{V}_{s+d_{1}f,g_{1}}(\,\b_{1};C(f)\,)\,
     \Phi^{V}_{s+d_{2}f,g_{2}}(\,\b_{2};C(pt)\,).
\end{align}
Therefore, if the limit set (\ref{limitset}) is contained in the
space (\ref{rel}) we have the sum formula (\ref{E:sumf}) from
(\ref{sumfm}), (\ref{left}) and (\ref{right}).

\smallskip In general, the limit set (\ref{limitset}) is not
contained in the space (\ref{rel}). In that case, there are maps
$f_{\l}$ in (\ref{cut-down1}) that converge to a limit $f_{0}$ as
$\l\to 0$ such that some components of $f_{0}$ map entirely into
$V$. The contribution of those maps in (\ref{cut-down1}) is call
the {\em contribution from the neck} and enters into the sum
formula (\ref{E:sumf}) as a correction term.  This correction term
can be computed by using the $S$-matrix (cf. section 12 of
\cite{ip3}). By the choice of the constraint $\b$ we have the
correction term
\begin{equation}\label{correction}
\sum\,     GW^{V}_{s+d_{1}f,g_{1}}(\,\b_{1};C(f)\,)\,
     \Phi^{V,U}_{d_{3}f,1}(\,C(pt),C(pt)\,)\,
     \Phi^{U}_{s+d_{2},g_{2}-1}(\,\b_{2};C(f)\,)
\end{equation}
where the sum is over all $d_{1}+d_{2}+d_{3}=d$. The correction
term (\ref{correction}) is zero by Lemma~\ref{rel-compu}\,e and
hence the proof is complete. \QED

\bigskip

Next, we use the sum formula (\ref{E:sumf}) to compute relative
invariants of $E(n)$.

\begin{lemma}\label{L:last}
Let $\gamma_{1}, \gamma_{2}$ be a basis of $H^{1}(\,E(0);Z\,)$.
\begin{enumerate}
  \item[(a)]
        $\Phi^{V}_{s+df,0}(\,\gamma_{1},\gamma_{1};C(f)\,)=1$
        if $d=0$ and 0 otherwise.
  \item[(b)]
        $\Phi^{V}_{s+df,1}(\,\gamma_{1},\gamma_{2};C(pt)\,)= 0$
  \item[(c)]
        $GW^{V}_{s+df,g}(\,pt^{g-1};C(pt)\,)\ =\ 0$
  \item[(d)]
        $GW^{V}_{s+df,g}(\,pt^{g};C(f)\,)\ =\
         GW^{\H}_{s+df,g}(\,pt^{g}\,)$
\end{enumerate}
\end{lemma}

\pf (a) follows from (i)
$\Phi_{s,0}(E(0))(\,\gamma_{1},\gamma_{2}\,)=1$\,(\,see Theorem 2
of \cite{ll}\,),  (ii) for $g=0$ relative invariants are same as
absolute invariants\,(\,Proposition 14.9 of \cite{ip3}\,), and
(iii) there is no rational curve representing $s+df$ with $d\ne 0$
on $E(0)=S^{2}\times T^{2}$ with a product complex structure

To prove (b), we will apply the sum formula (\,Theorem 12.4 of
\cite{ip3}\,) for the symplectic sum $E(0)=E(0)\#_{V}E(0)$. The
only difference between that sum formula and (\ref{E:sumf}) is the
degree of constraints. We also note that as in the proof of
Proposition~\ref{P:sumf} there is no contribution from the neck
for our case.

Split constraints $\gamma_{1}$ and $\gamma_{2}$ on one side and
one point constraint on the other side. Then we have
\begin{align}\label{last3}
\Phi_{s+df,1}(\,\gamma_{1},\gamma_{2},pt\,)\ &=
\sum_{d_{1}+d_{2}=d}
     \Phi^{V}_{s+d_{1}f,1}
     (\,\gamma_{1},\gamma_{2};C(pt)\,)\,
     \Phi^{V}_{s+d_{2},0}(\,pt;C(f)\,)
     \notag \\
     &+
     \sum_{d_{1}+d_{2}=d}
     \Phi^{V}_{s+d_{1}f,0}(\,\gamma_{1},\gamma_{2};C(f)\,)\,
     \Phi^{V}_{s+d_{2}f,1}(\,pt;C(pt)\,).
\end{align}
Using Lemma~\ref{rel-compu}\,c,d, and (a), we can simplify
(\ref{last3}) as
\begin{equation}\label{last33}
\Phi_{s+df,1}(\,\gamma_{1},\gamma_{2},pt\,)\ =\
\Phi^{V}_{s+df,0}(\,\gamma_{1},\gamma_{2};C(f)\,)\ +\
d\,\sigma(d).
\end{equation}
Now, (b) follows from (\ref{last33}) and
$\Phi_{s+df,1}(\,pt,\gamma_{1},\gamma_{2}\,)=d\,\sigma(d)$\, (see
Theorem 2 of \cite{ll}).

\medskip
Now, we use the sum formula (\ref{E:sumf}), (a) and (b) to show
(c) and (d). For the proof of (c), we split $g-1$ point
constraints on $E(n)$-side and the constraints $\gamma_{1}$ and
$\gamma_{2}$ on $E(0)$-side. Then we have
\begin{align}\label{last4}
&GW^{\H}_{s+df,g}(\,pt^{g-1},\pi_{\l}^{*}\gamma_{1},\pi_{\l}^{*}\gamma_{2}\,)
\notag\\
=& \sum_{d_{1}+d_{2}=d}
     GW^{V}_{s+d_{1}f,g}
     (\,pt^{g-1};C(pt)\,)\,
     \Phi^{V}_{s+d_{2},0}(\,\gamma_{1},\gamma_{2};C(f)\,)
     \notag \\
     +&
     \sum_{d_{1}+d_{2}=d}
     GW^{V}_{s+d_{1}f,g-1}(\,pt^{g-1};C(f)\,)\,
     \Phi^{V}_{s+d_{2}f,1}(\,\gamma_{1},\gamma_{2};C(pt)\,).
\end{align}
Since $E(n)$ is simply connected, the left-hand side of
(\ref{last4}) is zero. Therefore, (c) follows from (\ref{last4})
together with (a) and (b).

\smallskip
Lastly, we split $g$ point constraints on $E(n)$-side to obtain
  \begin{equation}\label{last5}
     GW^{\H}_{s+df,g}(\,pt^{g}\,)\ =
     \sum_{d_{1}+d_{2}=d}GW_{s+d_{1}f,g}^{V}(\,pt^{g};C(f)\,)\,
     \Phi_{s+d_{2}f,0}^{V}(\,C(pt)\,).
  \end{equation}
Now, (d) follows from (\ref{last5}) and Lemma~\ref{rel-compu}\,c.
\QED

\bigskip

Finally, we are ready to show the sum formulas (\ref{E:Sum1}) and
(\ref{E:Sum2}).

\begin{prop}{\bf (\,Sum Formulas\,)}
\label{Lastprop}
  \begin{enumerate}
     \item[(a)]
         $H(t)\  =\ 2\,F_{0}(t)\left(G(t)\ -\ \dfrac{1}{24} \right)$
     \item[(b)]
         $F_{g}(t)\ =\ F_{g-1}(t)\,t\,G^{\prime}(t)$
  \end{enumerate}
\end{prop}

\pf Choose a smooth fiber $F$ on $E(0)$-side of $Z_{0}$ and
consider the cut-down moduli space
\begin{equation*}
\M_{1,1}^{V}(\,d\,)\,\cap\,\tau(f^{*})\ =\ \{\ (h,\a)\in
\M_{1,1}^{V}(\,d\,)\ |\ ev(h)\in F\ \}.
\end{equation*}
The constraint $\tau(f^{*})$ lies only on $E(0)$-side. We can thus
apply the same argument as in the proof of
Proposition~\ref{P:sumf} to obtain
\begin{align}\label{E:sumf-a}
GW^{\H}_{s+df,1}(\,\tau(f^{*})\,)\ &= \sum_{d_{1}+d_{2}=d}
     GW^{V}_{s+d_{1}f,1}
     (\,C(pt)\,)\,
     \Phi^{V}_{s+d_{2},0}(\,\tau(f^{*});C(f)\,)
     \notag \\
     &+
     \sum_{d_{1}+d_{2}=d}
     GW^{V}_{s+d_{1}f,0}(\,C(f)\,)\,
     \Phi^{V}_{s+d_{2}f,1}(\,\tau(f^{*});C(pt)\,).
\end{align}
By Lemma~\ref{rel-compu}\,a,b and Lemma~\ref{L:last}\,d,
(\ref{E:sumf-a}) becomes
\begin{equation}\label{E:sumfaa}
GW^{\H}_{s+df,1}(\,\tau(f^{*})\,)\ =\ \sum_{d_{1}+d_{2}=d}
     2\,GW^{\H}_{s+d_{1}f,0}\,\sigma(d_{2}).
\end{equation}
Now, (a) follows from (\ref{E:sumfaa}) and the definition of
$F_{0}(t)$, $H(t)$, and $G(t)$.

\medskip
To prove (b) we split  $g-1$ points on $E(n)$-side and one point
on $E(0)$-side. Then by (\ref{E:sumf}) we have
\begin{align*}
GW^{\H}_{s+df,g}(\,pt^{g}\,)\ &= \sum_{d_{1}+d_{2}=d}
     GW^{V}_{s+d_{1}f,g}
     (\,pt^{g-1};C(pt)\,)\,
     \Phi^{V}_{s+d_{2},0}(\,pt;C(f)\,)
     \notag \\
     &+
     \sum_{d_{1}+d_{2}=d}
     GW^{V}_{s+d_{1}f,g-1}(\,pt^{g-1};C(f)\,)\,
     \Phi^{V}_{s+d_{2}f,1}(\,pt;C(pt)\,).
\end{align*}
By Lemma~\ref{rel-compu}\,d and Lemma~\ref{L:last}\,c,d, this
becomes
\begin{equation}\label{E:sumfb}
GW^{\H}_{s+df,g}(\,pt^{g}\,)\ =\ \sum_{d_{1}+d_{2}=d}
     GW^{\H}_{s+d_{1}f,g-1}(\,pt^{g-1}\,)\,d_{2}\sigma(d_{2}).
\end{equation}
Together with the definition of $F_{g}(t)$ and $G(t)$,
(\ref{E:sumfb}) implies (b). \QED

\vskip 1cm

%%%%%%%%%%%%%%%%%%%%%%%%%%%%%%%%%%%%%%%%%%%%%%%%%%%%%%%%%%%%%%%%%%%%%%%%%%%%%%%
%\pagebreak


\begin{thebibliography}{}


%\bibitem[A]{a} N. Aronszajn, {\em A unique continuation theorem for
%solutions of elliptic partial differential
%equations or inequalities of the second order}, J. Math Pures Appl. {\bf
%9}, (1957), 235-249.


\bibitem[AC]{ac} E. Arbarello and M. Cornalba, {\em Calculating cohomology
groups of moduli spaces of curves via algebraic geometry},
Inst. Hautes \'{E}tudes Sci. Publ. Math. No. 88 (1998), 97-127.



%\bibitem[B]{b} L. Besse, {\em Einstein manifolds},
%                  Springer-Verlag, Berlin Heidelberg, 1987.

\bibitem[BL1]{bl1} J. Bryan and N. C. Leung,
                      {\em The enumerative geometry of K3 surfaces
                           and modular forms}, J. Amer. Math. Soc.
                      {\bf 13} (2000), 371-410.

\bibitem[BL2]{bl2} J. Bryan and N. C. Leung, {\em Genenerating functions
for the number of curves on abelian surfaces},
Duke Math. J. {\bf 99} (1999), no. 2, 311-328.

\bibitem[BL3]{bl3} J. Bryan and N. C. Leung, {\em Counting curves on
irrational surfaces}, Surveys in differential geometry:
differential geometry inspired by string theory, 313-339,
Surv. Diff. Geom., 5, Int. Press, Boston, MA, 1999.

\bibitem[BL4]{bl4} J. Bryan and N. C. Leung, private conversation.



%\bibitem[BF]{bf} K. Behrend and B. Fantechi, In Preparation.

%\bibitem[BPV]{bpv} W. Barth, C. Peters, and A. Van de Ven,
%                      {\em Compact complex surfaces},
%                       Springer-Verlag, Berlin Heidelberg, 1984.

%\bibitem[CH]{ch} L. Caporaso and J. Harris, {\em Counting plane curves in
%any genus}, Invent. Math. {\bf 131} (1998), 345-392.


%\bibitem[D]{d} S. K. Donaldson,
%{\em Yang-Mills invariants of 4-manifolds },
%In S.K. Donaldson and C.B. Thomas, editors,
%Geometry of low dimensional manifolds : Gauge Theory and Algebraic
%surfaces,
%number 150 in London Math. Soc. Lecture Note Series,
%Cambrige University Press, 1989.


%\bibitem[E]{e} Y. Eliashberg, {\em Invariants in contact topology},
%       Proceedings of the International
%Congress of Mathematicians, Vol. II (Berlin, 1998). Doc. Math. 1998,
%Extra Vol. II, 327-338.


%\bibitem[FM]{fm} R. Friedman, J.W. Morgan, {\em Smooth four-manifolds
%       and complex surfaces},
%       Springer-Verlag, Berlin Heidelberg, 1994.




%\bibitem[FO]{fo} K. Fukaya, K. Ono, {\em Arnold conjecture and
%Gromov-Witten invariants},  Topology{\bf 38} (1999),
%933-1048.


\bibitem[G]{g} E. Getzler,
{\em Topological recursion relations in genus 2},
In "Integrable systems and algebraic geometry(Kobe/Kyoto, 1997)."
World Sci. Publishing, River Edge, NJ, 198, pp 73-106.



%\bibitem[G]{g} L. G\"{o}ttsche, {\em A conjectural generating
%function for numbers of curves on surfaces}, preprint,
%alg-geom/9711012


%\bibitem[G]{g} M. Gromov,  {\em Pseudo holomorphic curves in symplectic
%manifolds}, Invent. Math. {\bf 82} (1985), 307-347.


%\bibitem[GH]{gh} P. Griffiths and J. Harris, {\em Principles of algebraic
%geometry}, J. Wiley, New York, 1978.

%\bibitem[GT]{gt} D. Gilbarg and N. Trudinger, {\em Elliptic partial
%differential equations of second order},
%Springer-Verlag, Berlin Heidelberg, 1983.



\bibitem[IP1]{ip1} E. Ionel and T.  Parker, {\em Gromov Invariants
and Symplectic Maps}, Math. Annalen, {\bf 314}, 127-158 (1999).


\bibitem[IP2]{ip2} E. Ionel and T.  Parker,  {\em Relative
Gromov-Witten Invariants}, Ann. Math. 157 (2003), 45-96.


\bibitem[IP3]{ip3} E. Ionel and T.  Parker,  {\em The Symplectic Sum
Formula for Gromov-Witten Invariants}, preprint, math. SG/0010217.


\bibitem[IS]{is} S. Ivashkovich and V. Shevchishin,
{\em Gromov compactness theorem for $J$-complex curves with boundary},
Internat. Math. Res. Notices {\bf 2000}, no. 22, 1167-1206.


\bibitem[KM1]{km1} M. Kontsevich, Y.I. Manin,
{\em Relations between the correlators of the topological
     sigma model coupled to gravity},
Commun. Math. Phys. 196 (1998), 385-398.



\bibitem[KM2]{km2} P. Kronheimer and T. Mrowka,
{\em Embedded surfaces and the
structure of
Donaldson's polynomial invariants}, J. Differential Geom.
{\bf 41} (1995), 573-734.



%\bibitem[KP]{kp} S. Kleiman and R. Piene,
%{\em Enumerating singular curves on surfaces},
%Algebraic Geometry : Hirzebruch 70 (Warsaw, 1998), 209-238,
%Contemp. Math., 241, Amer. Math. Soc., Providence, RI, 1999.




\bibitem[L]{l} J. Lee, {\em Family Gromov-Witten Invariants for
K\"{a}hler Surfaces}, preprint, math. SG/0209402.

\bibitem[LL]{ll} T.J. Li, A. Liu,
         {\em Counting curves on elliptic ruled surface}, Topology and
         Its Applications, 124 (2002), 347-353.


\bibitem[Lo]{lo} E. Looijenga, {\em Smooth Deligne-Mumford
compactifications by means of Prym level structures}, J. Alg. Geom. 3 (1994),
no. 2, 283-29



%\bibitem[LR]{lr} A.-M. Li, Y. Ruan, {\em Symplectic surgery and Gromov-Witten
%invariants of Calabi-Yau 3-folds I}, preprint, alg-geom/9803036.



\bibitem[LT]{lt} J. Li and G. Tian,  {\em Virtual moduli cycles and
Gromov-Witten invariants of general symplectic manifolds}, Topics in
symplectic $4$-manifolds (Irvine, CA, 1996), 47--83, First Int.
Press Lect. Ser., I, International Press, Cambridge, MA, 1998.





%\bibitem[L]{lock} R. Lockhart,  {\em Fredholm, Hodge and Liouville
%theorems on non-compact manifolds}, Trans. Amer. Math. Soc. , {\bf 301}
%         (1987) 1--35.



%\bibitem[M]{m} D. McDuff, {\em The local behavior of holomorphic curves
%in almost complex
%$4$-manifolds}, J. Differential Geom. {\bf 34} (1991), 143-164.


%\bibitem[MS]{ms} D. McDuff and D. Salamon, {$J$-holomorphic curves
%and quantum cohomology}, A.M.S.,
%Providence, R.I., 1994.

\bibitem[P]{p} T. Parker, {\em Compactified moduli spaces of
pseudo-holomorphic curves} Mirror symmetry, III (Montreal,
PQ, 1995), 77--113, AMS/IP Stud. Adv. Math., 10, Amer. Math. Soc.,
Providence, RI, 1999.

\bibitem[PW]{pw} T. Parker and J. Wolfson,
{\em Pseudo-holomorphic maps and bubble trees},
Jour. Geometric Analysis, {\bf 3} (1993) 63-98.

%\bibitem[R]{r} Y. Ruan, {\em Virtual neighborhood and pseudo-holomorphic
%curves}, Proceedings of 6th Gokova Geometry-Topology Conference,
%Turkish J. Math. 23 (1999), no. 1, 161--231.

\bibitem[RT1]{rt1} Y. Ruan and G. Tian, {\em A mathematical theory of quantum
cohomology}, J. Differential Geom. {\bf 42} (1995), 259-367.

\bibitem[RT2]{rt2}  Y. Ruan and G. Tian, {\em Higher genus symplectic
invariants and sigma models coupled with gravity},  Invent. Math. {\bf 130}
(1997),
455-516.

%\bibitem[S]{s} B. Siebert, {\em Symplectic Gromov-Witten invariants},
%New trends in algebraic geometry (Warwick, 1996),
%375-424, London Math. Soc. Lecture Note Ser., 264, Cambridge
%Univ. Press, Cambridge, 1999.

%\bibitem[SY]{sy}  Y.T Siu and S.T. Yau, {\em  Compact Kähler manifolds of
%positive
%bisectional curvature},   Invent. Math. {\bf 59} (1980), 189-204.


%\bibitem[T]{t} C. H. Taubes, {\em Counting
%pseudo-holomorphic curves in dimension four}, J. Diff. Geom. {\bf 44}
%(1996), 818--893.


%\bibitem[V]{v} I. Vainsencher, {\em Enumeration of n-fold tangent
%hyperplanes to a surface}, J. Alg. Geom. {\bf 4} (1995), 503-526


\bibitem[YZ]{yz} S.T. Yau and E. Zaslow, {\em BPS States, String Duality,
and Nodal Curves on K3}, Nuclear Phys. B {\bf 471} (1996), 503-512




\end{thebibliography}
\end{document}